\documentclass{article}
\usepackage{amsthm}
\usepackage{amssymb}
\usepackage{amsmath}
\usepackage{colonequals}
\usepackage{graphicx}
\newcommand{\GG}{{\mathcal G}}
\newcommand{\Ss}{{\mathcal S}}
\newcommand{\RR}{{\cal R}}

\newcommand{\QQ}{{\cal Q}}

\newcommand{\el}{\mbox{el}}

\newcommand{\surf}{\mbox{{\rm surf}}}
\newcommand{\gen}{\mbox{{\rm gen}}}
\newcommand{\claim}[2]{\begin{equation}\mbox{\parbox{\linewidth}{{\em #2}}}\label{#1}\end{equation}}
\newtheorem{theorem}{Theorem}[section]
\newtheorem{corollary}[theorem]{Corollary}
\newtheorem{lemma}[theorem]{Lemma}

\newcount\claimno
\def\newclaim#1#2{
   \global\advance\claimno by 1\relax
   \bigskip\noindent\rlap{\rm(\the\claimno)}\ignorespaces
   \global\expandafter\edef\csname CLAIMLABEL#1\endcsname{\the\claimno}\relax
   \hangindent=33pt\hskip30pt{\sl#2}\bigskip}
\def\mylabel#1{{\label{#1}}}
\def\junk#1{}

\begin{document}
\title{Three-coloring triangle-free graphs on surfaces IV.  
Bounding face sizes of $4$-critical graphs\thanks{18 January 2015.}}
\author{%
     Zden\v{e}k Dvo\v{r}\'ak\thanks{Computer Science Institute (CSI) of Charles University,
           Malostransk{\'e} n{\'a}m{\v e}st{\'\i} 25, 118 00 Prague, 
           Czech Republic. E-mail: {\tt rakdver@iuuk.mff.cuni.cz}.
	   Supported by the Center of Excellence -- Inst. for Theor. Comp. Sci., Prague, project P202/12/G061 of Czech Science Foundation.}
 \and
     Daniel Kr{\'a}l'\thanks{Faculty of Informatics,
     Masaryk University, Botanick\'a 68A, 602 00 Brno, Czech Republic, and
     Mathematics Institute, DIMAP and Department of Computer Science, University
     of Warwick, Coventry CV4 7AL, UK. E-mail: {\tt dkral@fi.muni.cz}}
 \and
        Robin Thomas\thanks{School of Mathematics, 
        Georgia Institute of Technology, Atlanta, GA 30332. 
        E-mail: {\tt thomas@math.gatech.edu}.
        Partially supported by NSF Grants No.~DMS-0739366 and DMS-1202640.}
}
\date{Apr 27, 2019}
\maketitle
\begin{abstract}
Let $G$ be a $4$-critical graph with $t$ triangles, embedded in a surface of genus $g$.
Let $c$ be the number of $4$-cycles in $G$ that do not bound a $2$-cell face.
We prove that $$\sum_{\mbox{$f$ face of $G$}} (|f|-4)\le \kappa(g+t+c-1)$$
for a fixed constant $\kappa$, thus generalizing and strengthening several known results.
As a corollary, we prove that every triangle-free graph $G$ embedded in a surface of genus $g$
contains a set of $O(g)$ vertices such that $G-X$ is $3$-colorable.
\end{abstract}

\section{Introduction}

This paper is a part of a series aimed at studying the $3$-colorability
of graphs on a fixed surface that are either triangle-free, or have their
triangles restricted in some way. Historically the first result in this direction is the following
classical theorem of Gr\"otzsch~\cite{grotzsch1959}.

\begin{theorem}\label{grotzsch}
Every triangle-free planar graph is $3$-colorable.
\end{theorem}

Thomassen~\cite{thom-torus,thomassen1995-34,ThoShortlist}
found three reasonably simple proofs of this claim. 
Recently, two of us, in joint work with Kawarabayashi~\cite{DvoKawTho}
were able to design a linear-time algorithm to $3$-color triangle-free
planar graphs, and as a by-product found perhaps a yet simpler proof
of Theorem~\ref{grotzsch}.  Kostochka and Yancey~\cite{koyan} gave a completely
different proof as a consequence of their results on critical graphs.
The statement of Theorem~\ref{grotzsch}
cannot be directly extended to any surface other than the sphere.
In fact, for every non-planar surface $\Sigma$ there are infinitely many
$4$-critical graphs that can be embedded in $\Sigma$ (a graph is $4$-critical if it is not $3$-colorable,
but all its proper subgraphs are $3$-colorable).

Gimbel and Thomassen~\cite{gimbel} obtained the following elegant characterization
of $4$-critical triangle-free graphs embedded in the projective plane.
A graph embedded in a surface is a {\em quadrangulation} if every face
is bounded by a cycle of length four.

\begin{theorem}\mylabel{thm:gimtho}
A triangle-free graph embedded in the projective plane is $3$-colorable if and only
if it has no subgraph isomorphic to a non-bipartite
quadrangulation of the projective plane.
\end{theorem}

However, on other surfaces the structure of triangle-free $4$-critical graphs appears more
complicated. Thus we aim not for a complete characterization, but instead for a quantitative bound that can be used in
applications---roughly, that each $4$-critical triangle-free graph embedded in a fixed surface
has only a bounded number of faces of length greater than $4$.  One difficulty that
needs to be dealt with is that the claim as stated is false: there exist $4$-critical plane
graphs with arbitrarily many faces of length $5$ and with exactly $4$ triangles,
see \cite{tw-klein} for a construction.  By replacing an edge in each of the triangles
by a suitably chosen triangle-free graph, one can obtain a $4$-critical triangle-free
graph embedded in a surface of bounded genus and with arbitrarily many $5$-faces.
However, the examples constructed in this way contain a large number of non-contractible $4$-cycles.
This turns out to be the case in general, as shown by the main result of this paper, the following.

\begin{theorem}\mylabel{thm:main}
There exists a constant $\kappa$ with the following property.  Let $G$ be a graph embedded in
a surface of Euler genus $g$.  Let $t$ be the number of triangles in $G$ and let $c$ be the number of
$4$-cycles in $G$ that do not bound a $2$-cell face.  If $G$ is $4$-critical, then
$$\sum_{\mbox{$f$ face of $G$}} (|f|-4)\le \kappa(g+t+c-1).$$
\end{theorem}

The first step towards this result was obtained in the previous paper of this series~\cite{trfree3},
where we showed that Theorem~\ref{thm:main} holds under the additional assumption that the girth of $G$
is at least five (this special case strengthens a theorem of Thomassen~\cite{thomassen-surf}).  The inclusion of triangles and non-facial $4$-cycles in the statement
is merely a technicality, as $(\le\!4)$-cycles can be eliminated by replacing their edges by suitable
graphs of girth $5$ at the expense of increasing the genus.  Hence, the main contribution of this paper
is dealing with $4$-faces.

Theorem~\ref{thm:main} is the cornerstone of this series.
We will use it in future papers to deduce the following results.
First, in~\cite{trfree5} (whose preliminary version appeared as~\cite{dkt}) we will rely on the fact that the bound in
Theorem~\ref{thm:main} is linear in $t$ to answer affirmatively a question
of Havel~\cite{conj-havel} whether every planar graph with triangles sufficiently 
far apart is $3$-colorable.
Second, we will use Theorem~\ref{thm:main} and another result about
$3$-coloring graphs with most faces bounded by $4$-cycles to design
a polynomial-time algorithm to test whether an input triangle-free
graph embedded in a fixed surface is $3$-colorable.
That settles a problem of Gimbel and Thomassen~\cite{gimbel} and
completes one of two missing steps in a research program
initiated by Thomassen~\cite{thomassen-surf}. The latter asks whether
for fixed integers $k$ and $q$ the $k$-colorability
of graphs of girth $q$ on a fixed surface can be tested in polynomial time.
(The other step concerns $4$-coloring graphs on a fixed surface,
and prospects for its resolution in the near future are not very 
bright at the moment.)
With additional effort we will be able to implement our algorithm
to run in linear time.
Third, we will show that every triangle-free graph with an embedding 
of large edge-width in an orientable surface is $3$-colorable.
That generalizes a theorem of Hutchinson~\cite{locplanq}, 
who proved it under the stronger assumption that every face is even-sided.

More immediately, we apply Theorem~\ref{thm:main} to prove another generalization of Gr\"otzsch's theorem.
Kawarabayashi and Thomassen~\cite{kawthorem} proved that there exists a function $f$ such that
every triangle-free graph $G$ embedded in a surface of Euler genus $g$ contains a set $X$ of at most $f(g)$ vertices
such that $G-X$ is $3$-colorable.  They prove the claim for a function $f(g)=O(g\log g)$, and believe that
using a significantly more involved argument, they can improve the bound to linear (which is the best possible,
as exemplified by a disjoint union of copies of the Gr\"otzsch graph).  Using the theory we develop,
it is easy to prove this claim.
\begin{theorem}\label{thm-linbound}
There exists a constant $\beta>0$ with the following property.  Every triangle-free graph $G$ embedded in
a surface of Euler genus $g$ contains a set $X$ of at most $\beta g$ vertices such that $G-X$ is $3$-colorable.
\end{theorem}

The proof of Theorem~\ref{thm:main} follows the method we developed in~\cite{trfree2,trfree3}:
we show that collapsing a $4$-face (a standard reduction used when dealing with embedded triangle-free graphs)
does not decrease a properly defined weight.  By induction (with the base case given by the result
of~\cite{trfree3} for graphs of girth $5$), this shows that the weight of every $4$-critical embedded graph is bounded.
As the contribution of each $(\ge\!5)$-face to the weight is positive and linear in its size, the bound of Theorem~\ref{thm:main} follows.
Several technical difficulties are hidden in this brief exposition, however we worked out the solutions for them in~\cite{trfree2,trfree3}.
The main obstacle not encountered before is the need to avoid creating new non-facial $4$-cycles in the reduction,
which we deal with in Lemma~\ref{lemma-obstacle}.

Definitions and auxiliary results from \cite{trfree2,trfree3} that we are going to need in the proof are introduced in Sections~\ref{sec-def} and \ref{sec-res}.
In Section~\ref{sec-maingen}, we prove a version of Theorem~\ref{thm:main} generalized to allow some of the vertices to be precolored.  Theorem~\ref{thm:main} follows as
a straightforward corollary as outlined in Section~\ref{sec-main}.  Theorem~\ref{thm-linbound} is proved in Section~\ref{sec-linbound}.

\section{Definitions}\label{sec-def}

All graphs in this paper are simple, with no loops or parallel edges.

A \emph{surface}
is a compact connected $2$-manifold with (possibly null) boundary.  Each component of the boundary
is homeomorphic to a circle, and we call it a \emph{cuff}.  For non-negative integers $a$, $b$ and $c$,
let $\Sigma(a,b,c)$ denote the surface obtained from the sphere by adding $a$ handles, $b$ crosscaps and
removing the interiors of $c$ pairwise disjoint closed discs.  A standard result in topology shows that
every surface is homeomorphic to $\Sigma(a,b,c)$ for some choice of $a$, $b$ and $c$.
Note that $\Sigma(0,0,0)$ is a sphere, $\Sigma(0,0,1)$ is a closed disk, $\Sigma(0,0,2)$ is a cylinder,
$\Sigma(1,0,0)$ is a torus, $\Sigma(0,1,0)$ is a projective plane and $\Sigma(0,2,0)$ is a Klein bottle.
The \emph{Euler genus} $g(\Sigma)$ of the surface $\Sigma=\Sigma(a,b,c)$ is defined as $2a+b$.
For a cuff $C$ of $\Sigma$, let $\widehat{C}$ denote an open disk with boundary $C$ such that $\widehat{C}$ is disjoint from $\Sigma$, and let $\Sigma+\widehat{C}$ be
the surface obtained by gluing $\Sigma$ and $\widehat{C}$ together, that is, by closing $C$ with a patch.
Let $\widehat{\Sigma}=\Sigma+\widehat{C_1}+\ldots+\widehat{C_c}$, where $C_1$, \ldots, $C_c$ are the cuffs of $\Sigma$,
be the surface without boundary obtained by patching all the cuffs.

Consider a graph $G$ embedded in the surface $\Sigma$; when useful, we identify $G$ with the topological
space consisting of the points corresponding to the vertices of $G$ and the simple curves corresponding
to the edges of $G$.  We say that the embedding is \emph{normal} if every cuff of $\Sigma$ is equal to a cycle in $G$,
and we call such a cycle a \emph{ring}.
Throughout the paper, all graphs are embedded normally.
A \emph{face} $f$ of $G$ is a maximal arcwise-connected subset of $\Sigma-G$.
We write $F(G)$ for the set of faces of $G$.
The boundary of a face is equal to a union of closed walks of $G$, which we call the \emph{boundary walks} of $f$.  

Consider a ring $R$.  If $R$ is a triangle and at most one vertex of $R$ has degree greater than two in $G$, we say that $R$
is a \emph{vertex-like ring}.  A ring with only vertices of degree two is \emph{isolated}.  For a vertex-like ring $R$ that is
not isolated, the \emph{main} vertex of $R$ is its vertex of degree greater than two.
A vertex $v$ of $G$ is a \emph{ring vertex} if $v$ is contained in a ring (i.e., $v$ is drawn in the boundary of $\Sigma$), 
and $v$ is \emph{internal} otherwise.  
A subgraph $H$ of $G$ is \emph{contractible} if there exists a closed disk $\Delta\subseteq \Sigma$ such that $H\subset \Delta$
(in particular, a cycle $K$ is contractible if there exists a closed disk $\Delta\subseteq \Sigma$ with boundary equal to $K$).
The subgraph $H$ \emph{surrounds a cuff $C$} if $H$ is not contractible in $\Sigma$, but it is contractible in $\Sigma+\widehat{C}$.
We say that $H$ \emph{surrounds a ring $R$} if $H$ surrounds the cuff incident with $R$.
We say that $H$ is \emph{essential} if it is not contractible and it does not surround any cuff.

Let $G$ be a graph embedded in a surface $\Sigma$, let the embedding be
normal, and let $\cal R$ be the set of rings
of this embedding.
In those circumstances we say that $G$ is a \emph{graph in $\Sigma$
with rings $\cal R$.}
Furthermore, some vertex-like rings are designated as \emph{weak vertex-like rings}.  

The \emph{length} $|R|$ of a ring $R$ is the number of vertices of $R$.\footnote{Let us comment on a bit of a notational discrepancy between this paper and the previous papers of this series.
In the previous papers, we have defined the length of vertex-like rings differently (the weak vertex-like rings
had length 0, the non-weak ones length 1, while in this paper both are defined to have length 3).  This was necessary for technical reasons that are not relevant
in this paper's arguments, and we opt for the simpler definition here to avoid the need to keep track of this special case.
The only place where this change of definitions matters is in Theorem~\ref{thm:treti} taken from~\cite{trfree3}.
With the current definition of the vertex-like ring length, we claim a slightly weaker (and thus also correct) statement than the one proved in~\cite{trfree3}.}
For a face $f$, by $|f|$ we mean the sum of the lengths of the boundary walks of $f$ (in particular, if an edge
appears twice in the boundary walks, it contributes $2$ to $|f|$).
For a set of rings $\RR$, let us define $\ell({\RR})=\sum_{R\in\RR} |R|$.

Let $G$ be a graph with rings $\cal R$.  Let $H=\bigcup {\cal R}$ and let $H'$ be a (not necessarily induced) subgraph of $G$ obtained from $H$
by, for each weak vertex-like ring $R$, removing the main vertex and one of the non-main vertices of $R$
(or by removing two vertices of $R$ if $R$ has no main vertex),
so that $H'$ intersects $R$ in exactly one non-main vertex.  A {\em precoloring} $\psi$ of $\cal R$ is a 
$3$-coloring of the graph $H'$.  
A precoloring of $\cal R$ {\em extends to a $3$-coloring of $G$}
if there exists a $3$-coloring $\varphi$ of $G$ such that $\varphi(v)=\psi(v)$ for every $v\in V(H')$.
The graph $G$ is {\em $\cal R$-critical} if $G\neq H$ and for every proper subgraph
$G'$ of $G$ that contains $H$, there exists a precoloring of ${\cal R}$ that extends
to a $3$-coloring of $G'$, but not to a $3$-coloring of $G$.  For a precoloring $\kappa$ of $\cal R$
the graph $G$ is {\em $\kappa$-critical} if $\kappa$ does not extend to a $3$-coloring of $G$,
but it extends to a $3$-coloring of every proper subgraph of $G$ that contains $\cal R$.

Let us remark that if $G$ is $\kappa$-critical for some $\kappa$, then it is $\cal R$-critical,
but the converse is not true (for example, consider a graph consisting of a single ring with two chords).
On the other hand, if $\kappa$ is a precoloring of the rings of $G$ that does not extend to a $3$-coloring of $G$, then
$G$ contains a (not necessarily unique) $\kappa$-critical subgraph.

Let $G$ be a graph in a surface $\Sigma$ with rings $\RR$.
A face is {\em open $2$-cell} if it is homeomorphic to an open disk.
A face is {\em closed $2$-cell} if it is open $2$-cell and
bounded by a cycle. A face $f$ is {\em semi-closed $2$-cell} if it is open $2$-cell,
and if a vertex $v$ appears more than once in the boundary walk of $f$, then it appears exactly twice,
$v$ is the main vertex of a vertex-like ring $R$ and the edges of $R$ form part of the boundary walk of $f$.
A face $f$ is {\em omnipresent} if it is not open $2$-cell and
each of its boundary walks is a cycle bounding
a closed disk $\Delta\subseteq \widehat{\Sigma}\setminus f$ containing exactly one ring.

\section{Auxiliary results}
\label{sec-res}

First, let us state several simple properties of critical graphs (proofs can be found in \cite{trfree2}).

\begin{lemma}
\mylabel{lem:i012}
Let $G$ be a graph in a surface $\Sigma$ with rings $\RR$. If $G$ is $\RR$-critical, then every internal vertex of $G$ has degree at least three.
\end{lemma}

\begin{lemma}
\mylabel{lem:crit3conn}
Let $G$ be a triangle-free graph in a surface $\Sigma$ with rings $\RR$.  Suppose that each component of $G$ is a planar graph
containing exactly one of the rings.  If $G$ is ${\RR}$-critical, then each component of $G$ is $2$-connected.
\end{lemma}

Throughout the rest of the paper, let $s:\{2,3,4,\ldots\}\to{\mathbb R}$ be the function defined by
$s(l)=0$ for $l\le 4$,
$s(5)= 4/4113$, 
$s(6)=72/4113$, 
$s(7)=540/4113$,
$s(8)=2184/4113$ and $s(l)=l-8$ for $l\ge9$.
Based on this function, we assign weights to the faces.
Let $G$ be a graph embedded in $\Sigma$ with rings $\RR$.
For a face $f$ of $G$, we define $w(f)=s(|f|)$ if $f$ is open $2$-cell
and $w(f)=|f|$ otherwise.  We define $w(G,{\cal R})$ as the sum of $w(f)$ over all faces
$f$ of $G$.

The main result of \cite{trfree2} bounds the weight of graphs embedded in the disk with one ring $R$
(we write $R$-critical instead of $\{R\}$-critical for briefness).

\begin{theorem}\mylabel{thm:diskgirth5}
Let $G$ be a graph of girth at least $5$ embedded in the disk with one ring $R$ of length $l\ge 5$.
If $G$ is $R$-critical, then $l\ge 8$ and $w(G,\{R\})\le s(l-3)+s(5)$.
\end{theorem}

This is extended to general surfaces in~\cite{trfree3}.
Let $\gen(g,t,t_0,t_1)$ be a function defined for non-negative integers $g$, $t$, $t_0$ and $t_1$
such that $t\ge t_0+t_1$ as
$$\gen(g,t,t_0,t_1)=120g+48t-4t_1-5t_0-120.$$
Let $\surf(g,t,t_0,t_1)$ be a function defined for non-negative integers $g$, $t$, $t_0$ and $t_1$
such that $t\ge t_0+t_1$ as
\begin{itemize}
\item $\surf(g,t,t_0,t_1)=\gen(g,t,t_0,t_1)+116-42t=8-4t_1-5t_0$ if $g=0$ and $t=t_0+t_1=2$,
\item $\surf(g,t,t_0,t_1)=\gen(g,t,t_0,t_1)+114-42t=6t-4t_1-5t_0-6$ if $g=0$, $t\le 2$ and $t_0+t_1<2$, and
\item $\surf(g,t,t_0,t_1)=\gen(g,t,t_0,t_1)$ otherwise.
\end{itemize}
A graph $G$ embedded in a surface $\Sigma$ with rings $\RR$ has {\em internal girth at least five} if every $(\le\!4)$-cycle in $G$
is equal to one of the rings.
Let $t_0(\RR)$ and $t_1(\RR)$ be the number of weak and non-weak vertex-like rings in $\RR$, respectively.

\begin{theorem}[{\cite[Theorem 6.2]{trfree3}}]\mylabel{thm:treti}
There exists a constant $\eta_0$ with the following property.  Let $G$ be a graph embedded in
a surface $\Sigma$ with rings $\RR$. If $G$ is $\RR$-critical and has internal girth at least five, then
$$w(G,{\RR})\le \eta_0\cdot\surf(g(\Sigma),|\RR|, t_0(\RR), t_1(\RR)) + \ell({\RR}).$$
\end{theorem}

Consider a graph $H$ embedded in a surface $\Pi$ with rings $\QQ$, and let $f$ be a face of $H$.
There exists a unique surface whose interior is homeomorphic to $f$, which we denote by $\Pi_f$.
Note that the cuffs of $\Pi_f$ correspond to the facial walks of $f$.

Let $G$ be a graph embedded in $\Sigma$ with rings $\RR$.
Let $J$ be a subgraph of $G$ and let $S$ be a subset of faces of $J\cup\bigcup\RR$ such that
$J$ is equal to the union of the boundaries of the faces in $S$.
We define $G[S]$ to be the subgraph of $G$ consisting of $J$
and all the vertices and edges drawn inside the faces of $S$.
Let $C_1,C_2,\ldots,C_k$ be the boundary walks of the faces in $S$.
We would like to view $G[S]$ as a graph with rings $C_1$, \ldots, $C_k$.
However, the $C_i$'s do not necessarily have to be disjoint, and they do not have to be
cycles or isolated vertices.
To overcome this difficulty, we proceed as follows:
Suppose that $S=\{f_1,\ldots, f_m\}$.  For $1\le i\le m$, let $\Sigma_i$ be
a surface with boundary $B_i$ such that $\Sigma_i\setminus B_i$ is homeomorphic to $f_i$ (i.e., $\Sigma_i$ is homeomorphic to $\Sigma_{f_i}$).
Let $\theta_i:\Sigma_i\setminus B_i\to f_i$ be a homeomorphism that extends to a continuous mapping 
$\theta_i:\Sigma_i\to\overline{f_i}$, where $\overline{f_i}$ denotes the closure of $f_i$.
Let $G_i$ be the inverse image of $G\cap \overline{f_i}$ under $\theta_i$.
Then $G_i$ is a graph normally embedded in $\Sigma_i$.  We say that the set of embedded graphs $\{G_i:1\le i\le m\}$ 
is a {\em $G$-expansion of $S$} (the subgraph $J$ is uniquely determined by $S$ as the union of the boundary walks
of the faces in $S$).
Note that there is a one-to-one correspondence between the boundary walks of the faces of $S$ and the rings of the graphs
in the $G$-expansion of $S$; however, each vertex of $J$ may be split to several copies.
For $1\le i\le m$, we let ${\RR}_i$ be the set of rings of $G_i$,
where each vertex-like ring $R$ is weak if and only if $R$ is also a weak vertex-like ring of $G$.
We say that the rings in ${\RR}_i$ are the {\em natural rings} of $G_i$.

We use the following basic property of critical graphs proved in~\cite{trfree3}.

\begin{lemma}
\mylabel{lem:surfcritical}
Let $G$ be a graph in a surface $\Sigma$ with rings $\RR$,
and assume that $G$ is $\RR$-critical.
Let $J$ be a subgraph of $G$ and $S$ be a subset of faces of $J\cup\bigcup\RR$
such that $J$ is the union of the boundary walks of the faces of $S$.
Let $G'$ be an element of
the $G$-expansion of $S$ and let ${\RR}'$ be its natural rings.
If $G'$ is not equal to the union of the rings in ${\RR}'$, then $G'$ is ${\RR'}$-critical.
\end{lemma}

A frequently used corollary of Lemma~\ref{lem:surfcritical} concerns the case that $J$ is a contractible cycle.

\begin{lemma}
\mylabel{lem:diskcritical}
Let $G$ be a graph in a surface $\Sigma$ with rings $\RR$,
and assume that $G$ is $\RR$-critical.
Let $C$ be a non-facial cycle in $G$ bounding an open disk $\Delta\subseteq \Sigma$,
and let $G'$ be the graph consisting of the vertices and edges of $G$ drawn in the closure of $\Delta$.
Then $G'$ may be regarded as a graph embedded in the disk with one ring $C$,
and as such it is $C$-critical.
\end{lemma}

Together with the following following result of~\cite{gimbel}, this implies that every contractible $(\le\!5)$-cycle in a
critical graph without contractible triangles bounds a face.

\begin{theorem}\mylabel{thm:six}
Let $G$ be a triangle-free graph embedded in a disk with one ring $R$ of length at most $6$.  If $G$ is $R$-critical,
then $|R|=6$ and all faces of $G$ have length $4$.  Furthermore, if $R=r_1r_2\ldots r_6$ is an induced cycle and $\psi$
is a precoloring of $R$ that does not extend to a $3$-coloring of $G$, then $\psi(r_1)=\psi(r_4)$, $\psi(r_2)=\psi(r_5)$,
and $\psi(r_4)=\psi(r_6)$.
\end{theorem}

Furthermore, criticality is also preserved when cutting the surface, which was again shown in~\cite{trfree3}.

\begin{lemma}\mylabel{lemma-crcon}
Let $G$ be a graph in a surface $\Sigma$ with rings $\RR$,
and assume that $G$ is $\RR$-critical.  Let $c$ be a simple closed curve in $\Sigma$
intersecting $G$ in a set $X$ of vertices.  Let $\Sigma_0$ be one of the surfaces obtained from $\Sigma$ by
cutting along $c$.  Let us split the vertices of $G$ along $c$, let $G'$ be the part of the resulting graph embedded in $\Sigma_0$,
let $X'$ be the set of vertices of $G'$ corresponding to the vertices of $X$ and let ${\RR}'\subseteq {\RR}$ be the the rings of $G$ that
are contained in $\Sigma_0$.
Let $\Delta$ be an open disk or a disjoint union of two open disks disjoint from $\Sigma_0$ such that the boundary of $\Delta$ is equal to the cuff(s) of $\Sigma_0$ corresponding to $c$.
Let $\Sigma'=\Sigma_0\cup \Delta$.
Let $Y$ consist of all vertices of $X'$ that are not incident with a cuff in $\Sigma'$.  For each $y\in Y$, choose an open disk $\Delta_y\subset \Delta$
such that the closures of the disks are pairwise disjoint and the boundary of $\Delta_y$ intersects $G'$ exactly in $y$.
Let $\Sigma''=\Sigma'\setminus \bigcup_{y\in Y}\Delta_y$.  For each $y\in Y$, add to $G'$ a triangle $R_y$ with $y\in V(R_y)$ tracing the boundary of $\Delta_y$,
and let ${\RR}''={\cal R'}\cup \{R_y:y\in Y\}$, where the rings $R_y$ are considered as non-weak vertex-like rings.
If $G'$ is not equal to the union of the rings in ${\RR}''$, then $G'$ is ${\RR}''$-critical.
\end{lemma}

In particular, if $G'$ is an component of a ${\RR}$-critical graph, ${\RR}'$ are the rings contained in $G'$ and
$G'$ is not equal to the union of ${\RR}'$, then $G'$ is ${\RR}'$-critical.

Furthermore, the following claim was proved in~\cite{trfree3} (in a stronger version allowing non-vertex-like rings of length $3$).

\begin{lemma}\label{lemma-critshort}
Let $G$ be an $\{R_1,R_2\}$-critical graph embedded in the cylinder, where each of $R_1$ and $R_2$ is a vertex-like ring.
If every cycle of length at most $4$ in $G$ is non-contractible, then $G$ consists of $R_1$, $R_2$ and an edge between them.
In particular, neither $R_1$ nor $R_2$ is weak.
\end{lemma}

We need to show that the same result holds even if we allow contractible $4$-cycles.
\begin{corollary}\label{cor-critshort}
Let $G$ be an $\{R_1,R_2\}$-critical graph embedded in the cylinder, where each of $R_1$ and $R_2$ is a vertex-like ring.
If $G$ contains no triangle distinct from the rings, then $G$ consists of $R_1$, $R_2$ and an edge between them.
In particular, neither $R_1$ nor $R_2$ is weak.
\end{corollary}
\begin{proof}
We proceed by induction, and thus we can assume that the claim holds for all graphs with less than $|V(G)|$ vertices.
Without loss of generality, we can assume that neither $R_1$ nor $R_2$ is weak.
Let $\varphi$ be a precoloring of $\{R_1,R_2\}$ that does not extend to a $3$-coloring of $G$.

Suppose first that the distance between $R_1$ and $R_2$ in $G$ is at least four.
By Lemmas~\ref{lemma-critshort} and \ref{lem:diskcritical} and by Theorem~\ref{thm:six}, we can assume that $G$ contains
a $4$-face $f=v_1v_2v_3v_4$.  For $i=1,2$, let $G_i$ denote the graph obtained from $G$ by identifying $v_i$ with $v_{i+2}$ and suppressing the arising parallel edges
(observe that $v_iv_{i+2}\not\in E(G)$, since $G$ does not contain triangles distinct from the rings, and thus $G_i$ does not
contain loops).  Since every $3$-coloring of $G_i$ extends to a $3$-coloring of $G$, but $\varphi$ does not extend to a $3$-coloring of $G$,
we conclude that $\varphi$ does not extend to $G_i$.  Since the distance between $R_1$ and $R_2$ in $G$ is at least $4$,
$G_i$ does not contain any edge between $R_1$ and $R_2$, and by the induction hypothesis it follows that $G_i$ contains a triangle.

Hence, for $i=1,2$, there exists a $5$-cycle $C_i$ in $G$ containing the path $v_iv_{i+1}v_{i+2}$;
that is, $v_i$ is joined to $v_{i+2}$ by a path of length exactly three in $G$.
However, since $G$ is embedded in the cylinder (a part of the plane), the existence of these paths implies that two adjacent vertices of $f$
have degree at least three and are contained in a triangle.  This is a contradiction, since the only triangles in $G$ are vertex-like rings.

Therefore, there exists a path $P$ of length at most three between the main vertices $r_1$ and $r_2$ of $R_1$ and $R_2$.
Let $P$ be the shortest such path.
Let $G_1$ be the graph obtained from $G$ by removing the non-main vertices of the rings, drawn in the sphere obtained
by patching the cuffs incident with $R_1$ and $R_2$.  Let $f$ be the face of the subgraph $P$ of $G_1$, and let
$G_2$ be the unique element of the $G_1$-expansion of $\{f\}$, with ring $C$ of length at most $6$.
Since $G$ is the shortest path between $r_1$ and $r_2$, and since $G_1$ is triangle-free, the cycle $C$ is induced.

Suppose for a contradiction that $G_2\neq C$; similarly to Lemma~\ref{lem:surfcritical}, we can observe that then $G_2$ is $C$-critical.
Since $G_2$ is triangle-free, Theorem~\ref{thm:six} implies this is only possible if $|C|=6$.
Let $P=r_1p_1p_2r_2$ and $C=r_1p'_1p'_2r_2p''_2p'_1$.  We can extend $\varphi$ to a proper $3$-coloring $\varphi'$ of $P$;
we let $\psi$ be a $3$-coloring of $C$ given by $\psi(r_1)=\varphi(r_1)$, $\psi(r_2)=\varphi(r_2)$,
$\psi(p'_1)=\psi(p''_1)=\varphi'(p_1)$ and $\psi(p'_2)=\psi(p''_2)=\varphi'(p_2)$.
By Theorem~\ref{thm:six}, $\psi$ extends to a $3$-coloring of $G_2$, which implies that $\varphi$ extends to a $3$-coloring of $G$,
which is a contradiction.

Therefore $G_2=C$, and thus $G_1=P$.  Since $G$ does not contain
internal vertices of degree two by Lemma~\ref{lem:i012}, it follows that $G$ consists of $R_1$, $R_2$, and the edge $r_1r_2$, as required.
\end{proof}

\section{The main result}\label{sec-main}

For technical reasons, we are going to prove the following generalization of Theorem~\ref{thm:treti} instead of Theorem~\ref{thm:main}.
We say that a graph $G$ embedded in a surface $\Sigma$ with rings $\RR$ is {\em internally triangle-free} if every triangle in $G$
forms a vertex-like ring.

\begin{theorem}\mylabel{thm-maingen}
There exists a constant $\eta$ with the following property.  Let $G$ be a graph embedded in
a surface $\Sigma$ with rings $\RR$. If $G$ is $\RR$-critical, internally triangle-free, and contains no non-contractible $4$-cycles, then
$$w(G,{\RR})\le \eta\cdot\surf(g(\Sigma),|\RR|, t_0(\RR), t_1(\RR)) + \ell({\RR}).$$
\end{theorem}

The following section is devoted to the proof of Theorem~\ref{thm-maingen}.  Here, we show how it implies Theorem~\ref{thm:main}.

\begin{proof}[Proof of Theorem~\ref{thm:main}]
Let $\eta$ be the constant of Theorem~\ref{thm-maingen}; we can assume that $\eta\ge 1/80$.  Let $\kappa=1600\eta/s(5)$.  Let $J\subseteq G$
be the union of all triangles and all $4$-cycles that do not bound a $2$-cell face.  Note that $|E(J)|\le 4c+3t$.
Let $S$ be the set of all faces of $J$.  Let $\{G_1,\ldots, G_k\}$ be the $G$-expansion of $S$, where
$G_i$ is embedded in a surface $\Sigma_i$ with rings $\RR_i$ for $1\le i\le k$.  Note that $G_i$ is either equal to its rings or
$\RR_i$-critical, and that $G_i$ has no vertex-like rings.  Let the graph $G'_i$ embedded in $\Sigma_i$ with rings $\RR'_i$ be obtained from $G_i$ by,
for each ring $R\in\RR_i$ of length at most $4$, subdividing an edge of $R$ by $5-|R|$ new vertices.  Observe that $G'_i$ is
either equal to its rings or $\RR'_i$-critical, triangle-free and contains no non-contractible $4$-cycles,
and that $\ell(\RR'_i)\le 5\ell(\RR_i)/3$.
By Theorem~\ref{thm-maingen}, we have
\begin{align*}
z_i&\colonequals \sum_{\mbox{$f$ face of $G'_i$}} (|f|-4)\\
&\le \frac{1}{s(5)} \sum_{\mbox{$f$ face of $G'_i$}} w(f)=\frac{1}{s(5)}w(G'_i,\RR'_i)\\
&\le \frac{120\eta(g(\Sigma_i)+|\RR'_i|-1)+\ell(\RR'_i)}{s(5)}\\
&\le \frac{120\eta(g(\Sigma_i)+\ell(\RR'_i)-1)}{s(5)},
\end{align*}
because $\eta\ge 1/80$ and $\ell(\RR'_i)\ge3|\RR'_i|$.
Note that $\sum_{i=1}^k g(\Sigma_i)\le g$ and $\sum_{i=1}^k \ell(\RR'_i)\le \frac{5}{3}\sum_{i=1}^k \ell(\RR_i)=\frac{10|E(J)|}{3}$.
We conclude that
\begin{align*}
\sum_{\mbox{$f$ face of $G$}} (|f|-4)&\le \sum_{i=1}^k z_i\\
&\le \frac{120\eta(g+10|E(J)|/3-1)}{s(5)}\\
&\le \frac{120\eta(g+40c/3+10t-1)}{s(5)}\\
&=\kappa\Bigl((g+c+t-1)+\frac{37}{40}-\frac{37}{40}g-\frac{1}{4}t\Bigr).
\end{align*}
Aksionov~\cite{aksenov} proved that every planar graph with at most three triangles is $3$-colorable; since $G$ is $4$-critical,
we conclude that either $g\ge 1$ or $t\ge 4$, and in either case $\tfrac{37}{40}-\tfrac{37}{40}g-\tfrac{1}{4}t\le 0$.  Therefore,
$$\sum_{\mbox{$f$ face of $G$}} (|f|-4)\le \kappa(g+c+t-1),$$
as desired.
\end{proof}

\section{The proof of Theorem~\ref{thm-maingen}}\label{sec-maingen}


First, we show a key lemma that enables us to reduce $4$-faces.  Let us give the definitions we need to state this lemma.

Let $\Pi$ be a surface with boundary and $c$ a simple curve intersecting the boundary of $\Pi$ exactly in its ends.
The topological space obtained from $\Pi$ by cutting along $c$ (i.e., removing $c$ and adding two new pieces of boundary
corresponding to $c$) is a union of at most two surfaces.  If surfaces $\Pi_1,\ldots, \Pi_k$ are obtained from $\Pi$ by repeating this construction,
we say that $\{\Pi_1,\ldots,\Pi_k\}$ is a \emph{fragmentation} of $\Pi$; let us remark that zero repetitions of the construction are allowed, i.e., $\{\Pi\}$ is also a fragmentation of $\Pi$.

Let $G$ and $G'$ be graphs embedded in a surface $\Sigma$ with rings $\RR$, and let $J$ be a subgraph of $G$ containing $\bigcup\RR$ as a proper subgraph.
Suppose that there exists a collection $\{S_f:f\in F(G')\}$ of pairwise disjoint non-empty sets $S_f$ of faces of $J$ such that
$F(J)=\bigcup_{f\in F(G')} S_f$ and
for every face $f\in F(G')$, the set $\{\Sigma_h:h\in S_f\}$ is (up to homeomorphism) a fragmentation of $\Sigma_f$.
In this situation, we say that $J$ together with this collection $\{S_f:f\in F(G')\}$ forms a \emph{cover of $G$ by faces of $G'$}.  We say that a set $S_f$ is \emph{non-trivial}
if either $|S_f|>1$, or $|S_f|=1$ and the unique face of $J$ in $S_f$ is not a face of $G$ (i.e., some vertices or edges of $G$ are drawn within it).
With respect to a fixed cover, we define the \emph{elasticity} $\el(f)$ of a face $f\in F(G')$ to be $\left(\sum_{h\in S_f} |h|\right)-|f|$.

\begin{figure}
\includegraphics[width=12cm]{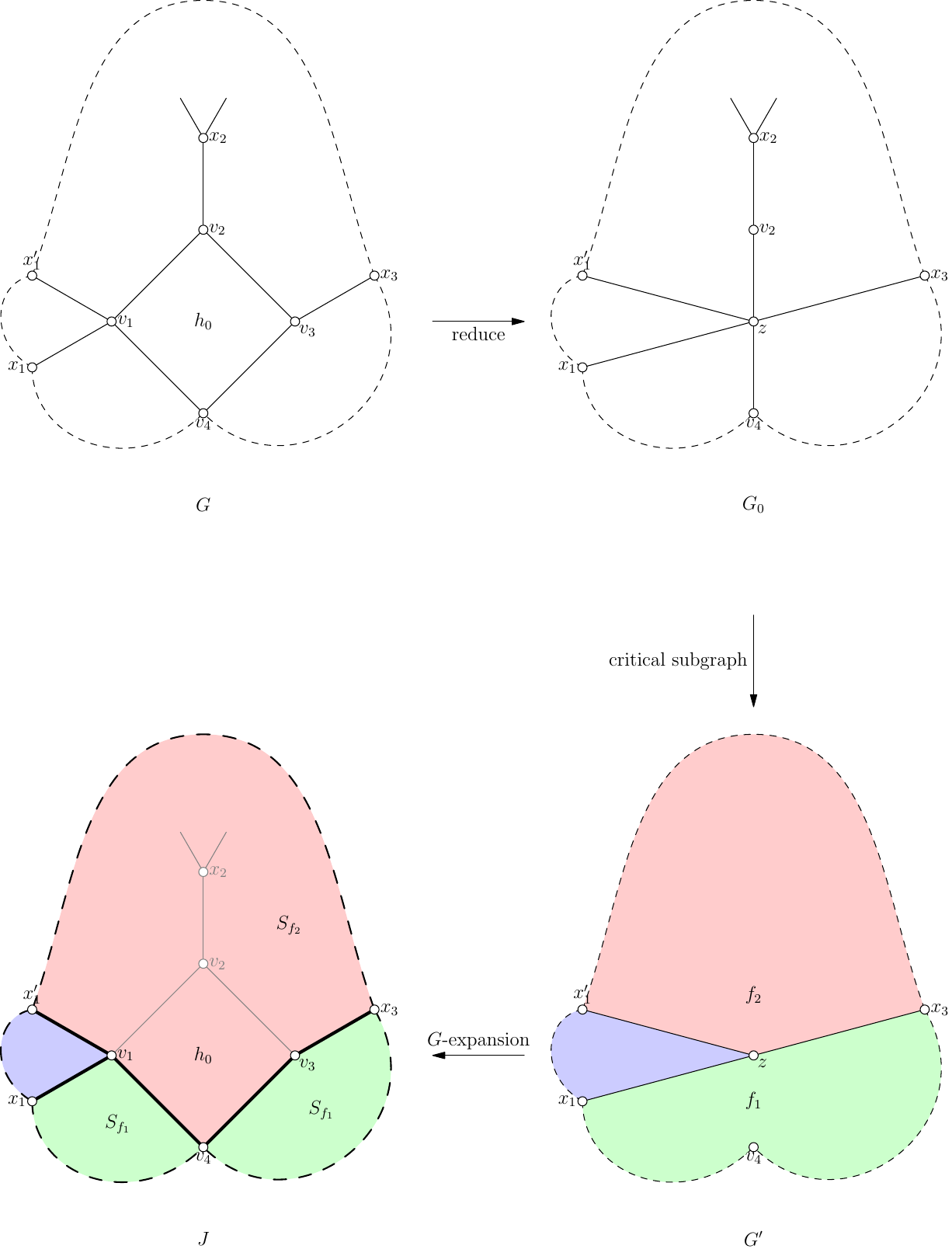}
\caption{The main idea.}\label{fig-mainstep}
\end{figure}

Before proceeding further, let us roughly explain what we aim to do, thus motivating the above definitions; see Figure~\ref{fig-mainstep} for an illustration.
We plan to prove Theorem~\ref{thm-maingen} by induction on the complexity of the surface $\Sigma$ and the size of $G$ (formalized by the definition of an ordering $\prec$, which we introduce later).
By Theorem~\ref{thm:treti}, we can assume $G$ contains
a closed 2-cell $4$-face $h_0$ bounded by a $4$-cycle $v_1v_2v_3v_4$.  Let $G_0$ be the graph obtained from $G$ by identifying $v_1$ with $v_3$ to a new vertex $z$.
Let us ignore for a moment a number of technical issues (triangles or non-contractible $4$-cycles could be created by this identification, both $v_1$ and $v_3$ could
be ring vertices) and suppose that $G_0$ is an internally triangle-free graph embedded in $\Sigma$ with the same rings $\RR$ as $G$ and without non-contractible $4$-cycles.  Let $\kappa$ be a precoloring of $\RR$
that does not extend to a $3$-coloring of $G$.  Since every $3$-coloring of $G_0$ extends to a $3$-coloring of $G$, we conclude that $\kappa$ does not extend to a $3$-coloring of $G_0$.
Therefore, although the graph $G_0$ does not have to be $\RR$-critical (e.g., in Figure~\ref{fig-mainstep}, the vertex $v_2$ has degree two in $G_0$, and thus it is irrelevant
for $3$-colorability), it must contain an $\RR$-critical subgraph $G'$.  By the induction hypothesis, we have $w(G',{\RR})\le \eta\cdot\surf(g(\Sigma),|\RR|, t_0(\RR), t_1(\RR)) + \ell({\RR})$,
and thus to prove Theorem~\ref{thm-maingen}, it suffices to argue that $w(G,\RR)\le w(G',\RR)$.

The graph $G'$ is not quite a subgraph of $G$; however, it can be transformed to a subgraph $J$ of $G$ by splitting the vertex $z$ back into $v_1$ and $v_3$ and adding the path $v_1v_4v_3$.
Each face $f$ of $G'$ then naturally corresponds to a set $S_f$ of faces of $J$---usually, this set has only one element; however, it can happen that a face of $G'$ is split
during the construction of $J$, see e.g. the face $f_1$ in Figure~\ref{fig-mainstep}, in which case $\{\Sigma_h:h\in S_f\}$ forms (up to homeomorphism) a fragmentation of $\Sigma_f$.
It is easy to see that $J$ together with the system $\{S_f:f\in F(G')\}$ forms a cover of $G$ by faces of $G'$.
Hence, to prove $w(G,\RR)\le w(G',\RR)$, we need to argue that for each $f\in F(G')$, the total weight of the faces of $G[S_f]$ is at most $w(f)$.

Let $J_f$ denote the union of the boundaries of the faces in $S_f$.  It is easy to see that since $G$ is $\RR$-critical, its subgraph $G[S_f]$ is critical with respect to $J_f$.
There is a caveat: we may not be able to directly embed $G[S_f]$ in a surface (or surfaces, if $|S_f|>1$) with rings corresponding to $J_f$; e.g., suppose that $\Sigma$ is
a torus, and $J_f$ is the union of two homotopic non-contractible cycles intersecting in one vertex $w$, and $S_f$ consists of the face $h$ of $J_f$ homeomorphic to the open cylinder.
Then the two boundary walks of $h$ intersect, but distinct rings of a graph must be vertex-disjoint.  However, we can easily overcome this difficulty by splitting the vertex
$w$ into two vertices, so that $J_f$ becomes a disjoint union of cycles and $G[S_f]$ can be naturally embedded in a cylinder with two rings (which is the surface
we below denote by $\Sigma_h$).  We formalized this splitting operation in Section~\ref{sec-res} as the $G$-expansion of $S_f$, and argued it preserves criticality in Lemma~\ref{lem:surfcritical}.

Nevertheless, ignoring this subtlety, since $G[S_f]$ is $J_f$-critical, we can apply the induction hypothesis to it (ignoring yet another technical issue: to justify this step, 
we need to argue that the surfaces $\Sigma_h$ for $h\in S_f$ are not more complex than $\Sigma$).  If the sum of the lengths of the rings of $G[S_f]$ is the same as $|f|$, this directly
gives the desired bound on the weights of faces in $G[S_f]$ in terms of $w(f)$.  However, this does not need to be the case for the faces touching $z$; in Figure~\ref{fig-mainstep},
the face in $S_{f_2}$ has length $|f_2|+2$ (in terms defined above, $f_2$ has elasticity $2$).  A slightly more involved accounting is needed for such faces of non-zero elasticity.

Let us now proceed with a formal argument.  As we indicated above, the described reduction is not valid for all choices of the $4$-face $h_0$.  In particular, it may fail when $h_0$
contains ring vertices.
Let $G$ be a graph embedded in a surface $\Sigma$ with rings $\RR$.
A $4$-face $h_0=v_1v_2v_3v_4$ is \emph{ring-bound} if (up to a relabelling of the vertices of $f$)
\begin{itemize}
\item $v_1$ is the main vertex of a vertex-like ring $R$ and there exists a cycle $C$ of length at most $6$ surrounding $R$ such that $v_2v_3v_4\subset C$; or,
\item $h_0$ is incident with vertices of two distinct rings; or,
\item $v_1$ and $v_3$ are ring vertices.
\end{itemize}

\begin{figure}
\begin{center}
\includegraphics[width=110mm]{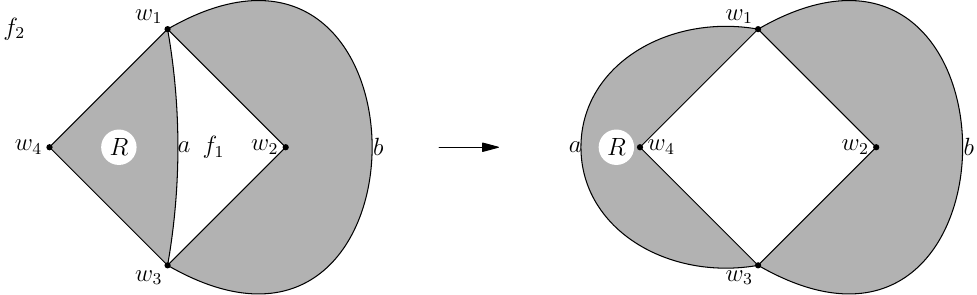}
\end{center}
\caption{Non-contractible $4$-cycle $w_1w_2w_3w_4$ flippable to a $4$-face, and the corresponding transformation.}\label{fig-flip}
\end{figure}

Another reason for the reduction to fail is that the identification of $v_1$ with $v_3$ creates a non-contractible $4$-cycle.  We cannot completely prevent this from
happening, but we can at least ensure that such a $4$-cycle has a quite special form.
A non-contractible $4$-cycle $C=w_1w_2w_3w_4$ is \emph{flippable to a $4$-face} if $C$ surrounds a ring $R$, $w_1w_2w_3$ is part of a boundary walk of a face $f_1$,
$w_1w_4w_3$ is part of a boundary walk of a face $f_2$, $C$ separates $f_1$ from $f_2$, and $w_1w_2w_3$ and $w_1w_4w_3$ are the only paths of length at most $2$ in $G$
between $w_1$ and $w_3$.  See the left part of Figure~\ref{fig-flip} for an illustration.

\begin{lemma}\label{lemma-obstacle}
Let $G$ be a graph embedded in a surface $\Sigma$ with rings $\RR$.  Suppose that $G$ is $\RR$-critical and internally triangle-free, contains no non-contractible $4$-cycles,
and every connected essential subgraph of $G$ has at least $13$ edges.
If $G$ has a $4$-face that is not ring-bound, then there exists an internally triangle-free $\RR$-critical graph $G'$ embedded in $\Sigma$ with rings $\RR$
such that $|E(G')|<|E(G)|$, and a subgraph $J\subseteq G$ containing $\bigcup\RR$ which together with a collection $\{S_f:f\in F(G')\}$ forms a cover of $G$ by faces of $G'$,
satisfying the following conditions.
\begin{itemize}
\item The graph $G'$ has at most one non-contractible $4$-cycle, and if it has one, it is flippable to a $4$-face.
\item The cover satisfies $\sum_{f\in F(G')} \el(f)\le 4$, and if a face $f\in F(G')$ is omnipresent or semi-closed $2$-cell,
then $\el(f)\in\{0,2\}$.
\item For every face $f\in F(G')$, if $f$ is semi-closed $2$-cell and $\el(f)=2$, then $S_f$ is non-trivial.
\end{itemize}
\end{lemma}
\begin{proof}
Let $h_0=v_1v_2v_3v_4$ be a $4$-face of $G$ that is not ring-bound.
In particular, $v_1$ and $v_3$ are not both ring vertices, and $v_2$ and $v_4$ are not both ring vertices;
by symmetry, we can assume that $v_3$ and $v_4$ are internal vertices.
By Lemma~\ref{lem:i012}, there exist edges $e_3$ and $e_4$ incident with $v_3$ and $v_4$, respectively,
and not incident with $h_0$.

As $G$ is internally triangle-free, $v_1$ and $v_3$ are non-adjacent
and $v_2$ and $v_4$ are non-adjacent.  For $i\in\{1,2\}$, we say that the vertex $v_i$ is \emph{problematic}
if either $v_i$ is the main vertex of a vertex-like ring, or
$G$ contains a path $P_i$ of length at most four joining $v_i$ with $v_{i+2}$ such
that the concatenation of the paths $P_i$ and $v_iv_{5-i}v_{i+2}$ is a non-contractible cycle $C_i$.
By switching the labels of vertices if necessary, we can assume that if $v_1$ is problematic, then $v_2$ is problematic as well.

Suppose that $v_1$ is the main vertex of a vertex-like ring $R$.  Then $v_1$ is problematic, and thus $v_2$ is problematic.  As $h_0$ is not ring-bound,
$v_2$ is an internal vertex, and thus $G$ contains a non-contractible cycle $C_2$ described above.  Since $h_0$ is not ring-bound,
we conclude $C_2$ does not surround
$v_1$.  However, then $C_2$ together with $R$ and the edge $v_1v_2$ forms a connected essential subgraph of $G$ with at most $11$ edges,
which is a contradiction.
We conclude that $v_1$ is not the main vertex of a vertex-like ring; and by symmetry, if both $v_1$ and $v_2$ are problematic,
then $v_2$ is not the main vertex of a vertex-like ring.

Let $G_0$ be a multigraph obtained from $G$ by deleting edges $v_2v_3$ and $v_3v_4$ and identifying $v_1$ with $v_3$ to a new vertex $z$
(parallel edges may arise if $v_1$ and $v_3$ have common neighbors other than $v_2$ and $v_4$).
Since the vertex $v_3$ is internal, $G_0$ is embedded in $\Sigma$ with rings $\RR$.  Furthermore, note that each edge $e\in E(G_0)$ corresponds
to a unique edge $\overline{e}$ of $G$.  Every $3$-coloring of $G_0$ corresponds to a $3$-coloring of $G$ in which $v_1$ and $v_3$ have the same color as $z$.
Let $\psi$ be a precoloring of $\RR$ that does not extend to $G$, and observe that $\psi$ also does not extend to $G_0$.
Therefore, $G_0$ has a $\psi$-critical subgraph $G'$ (the criticality implies $G'$ does not have parallel edges).

If there exists a path $xzy$ in $G_0$ contained in the boundary of a face such that $\overline{xz}$ is incident with $v_1$ in $G$
and $\overline{yz}$ is incident with $v_3$ in $G$, then let $A=\{v_1v_4,v_3v_4\}$, otherwise let $A=\emptyset$.
Let $J$ be the subgraph of $G$ with edge set $\{\overline{e}:e\in E(G')\}\cup A$ and with the vertex set consisting of all vertices
incident with these edges (informally, we split the vertex $z$ back into $v_1$ and $v_3$, and if that breaks the boundary of any face,
we fix it by adding the path $v_1v_4v_3$).  Each face $f$ of $G'$ naturally corresponds to either one or two faces of $J$
(the latter happens when $v_4\in V(G')$, $z$ is not adjacent to $v_4$ in $G'$, and $A\neq\emptyset$, see the face $f_1$ in
Figure~\ref{fig-mainstep} for an illustration); let $S_f$ denote the set of faces of $J$ corresponding to $f$.
The construction implies that $\{\Sigma_h:h\in S_f\}$ is (up to homeomorphism) a fragmentation of $\Sigma_f$.
Furthermore, every face of $J$ is contained in $S_f$ for a unique face $f\in F(G')$.
Therefore, $J$ together with the collection $\{S_f:f\in F(G')\}$ forms a cover of $G$ by faces of $G'$.

Observe that a face $f\in F(G')$ has non-zero elasticity if and only if 
if a facial walk of $f$ contains two consecutive edges $e_1$ and $e_2$ incident with $z$ such that $\overline{e_1}$ is incident with $v_1$ in $G$,
$\overline{e_2}$ is incident with $v_3$ in $G$, and $\overline{e_1}\neq v_1v_4$; and each such pair of edges contributes $2$ to the elasticity
of $f$.  There are at most two such pairs of edges in the cyclic order around $z$, and thus the sum of the elasticities of faces of $G'$ is at most $4$.
Furthermore, if both of these pairs of edges are incident with the same face $f$, then $z$ is incident twice with $f$;
since $z$ is not the main vertex of a vertex-like ring of $G'$ (as neither $v_1$ nor $v_3$ is the main vertex of a vertex-like ring
in $G$, and $G$ has no non-vertex-like rings of length $3$),
it follows that the face $f$ of $G'$ is neither semi-closed $2$-cell nor omnipresent.  We conclude that if a face $f\in F(G')$ is semi-closed $2$-cell
or omnipresent, then its elasticity is at most two.

Suppose that $f$ is a semi-closed $2$-cell face of $G'$, the elasticity of $f$ is not $0$, and $|S_f|=1$ (let $h$ be the unique face in $S_f$).
Then the boundary of $h$ contains the path $v_1v_4v_3$.  Note that $v_2v_3\not\in E(J)$, and thus if $h_0$ is a subset of $h$,
then an edge of $G$ is drawn within $h$.  If $h_0$ is not a subset of $h$, then the edge $e_4$ of $G$
is drawn within $h$.  We conclude that for every semi-closed $2$-cell face $f$ of $G'$ with $\el(f)>0$, the set $S_f$ is non-trivial.

To finish the proof of Lemma~\ref{lemma-obstacle}, it suffices to argue that $G'$ is internally triangle-free,
has at most one non-contractible $4$-cycle, and if it has one, it is flippable to a $4$-face.
Hence, suppose that $G'$ contains either a triangle which is not a vertex-like ring, or a non-contractible $4$-cycle.
Then $G$ contains a path $P_1$ of length $3$ or $4$ joining $v_1$ with $v_3$ such that $P_1$ together
with $v_1v_4v_3$ forms a cycle $C_1$ such that if $|C_1|=6$, then $C_1$ is non-contractible.  Let us first consider the case that $|C_1|=5$ and
$C_1$ is contractible, and let $\Delta\subset\Sigma$ be a closed disk bounded by $C_1$.  Since $v_1v_4v_3$ is a subpath of $C$,
either the edge $e_4$ or the $4$-face $h_0$ of $G$ is contained in the interior of $\Delta$.  Consequently, Lemma~\ref{lem:diskcritical} implies that the subgraph of $G$ drawn
in $\Delta$ is $C_1$-critical.  However, as $|C_1|=5$, this contradicts Theorem~\ref{thm:six}.

We conclude that $C_1$ is non-contractible, and thus $v_1$ is problematic.  By the choice of the labels of vertices incident with $h_0$,
it follows that $v_2$ is problematic as well.  We argued before that if both $v_1$ and $v_2$ are problematic, then $v_2$ is not the main
vertex of a vertex-like ring; and thus there exists a path $P_2$ joining $v_2$ with $v_4$ which combines
with the path $v_2v_3v_4$ to a non-contractible $(\le\!6)$-cycle $C_2$.

Note that the connected subgraph $H$ of $G$ formed by the cycle $v_1v_2v_3v_4$ and the paths $P_1$ and $P_2$ has at most $12$ edges,
and thus it is not essential.  It follows that there exists a cuff $C$ (incident with a ring $R$)
and a closed disk $\Lambda\subseteq \Sigma+\widehat{C}$ with $H$ contained in $\Lambda$.
Observe that $H$ does not contain a vertex-like ring, and since $G$ is internally triangle-free, it follows that $H$ is triangle-free.
Extend the disk $\Lambda$ to a plane and consider the resulting drawing of $H$ in the plane; for $1\le i\le 4$, let $h_i$ be the face of $H$ in this plane drawing
incident with $v_iv_{i+1}$ (where $v_5=v_1$) and distinct from $h_0$.  The faces $h_1$, \ldots, $h_4$ are pairwise distinct,
as any two of them are separated by $C_1$ or $C_2$.
Since $H$ is triangle-free, all faces of $H$ have length at least four.  Furthermore, observe that
the boundary walks of two faces of $H$ (the face containing $\widehat{C}$ and the unbounded face) surround $R$,
and since $G$ does not contain non-contractible $4$-cycles, it follows that these two faces have length at least $5$.
Therefore,
\claim{cl-facin}{$$22\le \sum_{h\in F(H)} |h|=8+2|E(P_1)\cup E(P_2)|\le 24.$$}
In particular, $H$ has no face distinct from $h_0$, \ldots, $h_4$.
By symmetry, we can assume that $\widehat{C}$ is contained in $h_3$, and thus $|h_3|\ge 5$.
The cycles $C_1$ and $C_2$ are non-contractible, and thus they surround $R$;
since $h_1$ is the only face of $H$ distinct from $h_0$ separated from $h_3$ by both $C_1$ and $C_2$,
it follows that $h_1$ is the unbounded face of $H$.  Hence, the boundary walk of $h_1$ surrounds $R$ and $|h_1|\ge 5$.
Note that $|h_1|+|h_3|\le |E(P_1)|+|E(P_2)|+2\le 10$, and thus $|h_1|=|h_3|=5$ and $|E(P_1)|=|E(P_2)|=4$.
As the subdisk of $\Lambda$ bounded by $C_2$ consists of the closure of the faces $h_2$ and $h_3$
and $|C_2|=6$, it follows that $|h_2|$ and $|h_3|$ have the same parity, i.e., $|h_2|$ is odd.
Similarly, $|h_4|$ is odd.  We conclude that $|h_2|,|h_4|\ge 5$, and by (\ref{cl-facin}),
it follows that $|h_2|=|h_4|=5$ and the paths $P_1$ and $P_2$ are edge-disjoint.  That is,
the $P_1$ and $P_2$ intersect in a single vertex $w_1$.  Observe that since $|h_2|=|h_3|=|h_4|$,
$w_1$ is the middle vertex of $P_1$ and $P_2$.

Since $C_1$ and $C_2$ surround $R$, the cycles bounding $h_2$ and $h_4$ are contractible, and thus by Lemma~\ref{lem:diskcritical} and Theorem~\ref{thm:six}, $h_2$ and $h_4$ are faces of $G$.
Let $Q_1=v_1w_2w_1$ be the subpath of $P_1$ between $v_1$ and $w_1$.  If $G$ contained another path $Q'_1$ of length at most two between $v_1$ and $w_1$,
then the cycle consisting of $Q_1$ and $Q'_1$ would be a non-ring triangle, or a non-contractible $4$-cycle, or a $4$-face showing that $w_2$ is an internal
vertex of degree two, which is a contradiction.  Therefore, $Q_1$ is the only path of length at most two between $v_1$ and $w_1$ in $G$.
Similarly, the subpath $v_3w_4w_1$ of $P_2$ is the only path of length at most two between $v_3$ and $w_1$ in $G$.

We conclude that that the graph $G'$ is internally triangle-free and contains (at most) one non-contractible $4$-cycle
$w_1w_2w_3w_4$, where $w_3=z$.  Because of the faces $h_2$ and $h_4$ of $G$, the paths $w_1w_2w_3$
and $w_1w_4w_3$ are parts of boundaries of faces of $G'$ separated by $w_1w_2w_3w_4$.  Furthermore, $w_1w_2w_3w_4$ surrounds the ring $R$.
We conclude that $w_1w_2w_3w_4$ is flippable to a $4$-face.
\end{proof}

As the next step, let us prove a stronger variant of Theorem~\ref{thm-maingen} in the case $\Sigma$ is the disk.
Let $\GG_{r,k}$ denote the class of all $R$-critical graphs of girth at least $r$ embedded in the disk with one ring $R$ of length $k$.
For a graph $G\in\GG_{4,k}$ let $S(G)$ denote the multiset of lengths of the $(\ge\!5)$-faces of $G$.  
Let $\Ss_{r,k}=\{S(G):G\in\GG_{r,k}\}$.  Note that by Theorem~\ref{thm:diskgirth5}, the set $\Ss_{5,k}$ is empty for $5\le k\le 7$
and finite for every $k\ge 8$,
the maximum of each element of $\Ss_{5,k}$ is at most $k-3$, and if the maximum of $S\in\Ss_{5,k}$ is equal to $k-3$,
then $S=\{5,k-3\}$.
Furthermore, $\Ss_{4,4}=\Ss_{4,5}=\emptyset$ and $\Ss_{4,6}=\{\emptyset\}$ by Theorem~\ref{thm:six}.

Let $S_1$ and $S_2$ be multisets of integers.  We say that $S_2$ is a \emph{one-step refinement} of $S_1$
if there exists $k\in S_1$ and a multiset $Z\in \Ss_{4,k}\cup \Ss_{4,k+2}$ such that $S_2=(S_1\setminus \{k\})\cup Z$.
We say that $S_2$ is a \emph{refinement} of $S_1$ if it can be obtained from $S_2$ by a (possibly empty)
sequence of one-step refinements.

\begin{lemma}\label{lemma-diskcase}
For every $k\ge 7$, each element of $\Ss_{4,k}$ other than $\{k-2\}$ is a refinement of an element of $\Ss_{4,k-2}\cup \Ss_{5,k}$.
In particular, if $S\in\Ss_{4,k}$, then the maximum of $S$ is at most $k-2$, and if the maximum is equal to $k-2$,
then $S=\{k-2\}$.
\end{lemma}
\begin{proof}
Let us re-phrase the claim.
\claim{cl-rep}{For every graph $G\in \GG_{4,k}$, either $S(G)=\{k-2\}$ or $S(G)$ is a refinement of an element of $\Ss_{4,k-2}\cup \Ss_{5,k}$.}
We prove this reformulation by induction on the number of edges of $G$; that is, we assume that (\ref{cl-rep})
holds for every graph $G'\in\GG_{4,k}$ with $|E(G')|<|E(G)|$.

Let us remark that this implies that the maximum of each element of $\Ss_{4,k}$ distinct from $\{k-2\}$
is at most $k-3$, by induction on $k$. Indeed, assuming this claim holds for every $k'<k$, the maximum of each element $\Ss_{4,k-2}\cup \Ss_{5,k}$
is at most $k-3$ by Theorem~\ref{thm:diskgirth5}, and thus it suffices to argue that this property is preserved by one-step refinements.
This is the case, since in the one-step refinement operation, we are replacing an element $c\le k-3$ by elements
of $\Ss_{4,c}\cup \Ss_{4,c+2}$, which are all smaller or equal to $(c+2)-2\le k-3$ by the assumption.

Let us now proceed with the proof of (\ref{cl-rep}).
By Theorem~\ref{thm:six}, (\ref{cl-rep}) holds if $k\le 6$, and thus we can assume that $k\ge 7$.
Consider a graph $G\in\GG_{4,k}$.
Note that every $4$-cycle in $G$ bounds a face by Lemma~\ref{lem:diskcritical} and Theorem~\ref{thm:six}.
If $G$ has no $4$-faces, then $S(G)\in \Ss_{5,k}$ and (\ref{cl-rep}) holds.  Hence, suppose $G$ contains a $4$-face $h_0=v_1v_2v_3v_4$.

Let us first consider the case that $h_0$ is ring-bound. Since $G$ has only one ring and this ring is not vertex-like,
we can assume both $v_1$ and $v_3$ belong to the ring $R$ of $G$. Consequently $G$ contains a path $P$ of length $1$ or $2$
intersecting $R$ exactly in the endpoints of $P$.
Let $C_1$ and $C_2$ be the cycles of $R+P$ distinct from $R$, and let $G_1$ and $G_2$ be the subgraphs of $G$ drawn in the closed
disks bounded by $C_1$ and $C_2$, respectively.  We have $|C_1|+|C_2|=k+2|E(P)|$.  Suppose that say $|C_1|=4$; hence,
$C_1$ bounds a face by Lemma~\ref{lem:diskcritical} and Theorem~\ref{thm:six}, and $S(G)=S(G_2)$.
If $C_2$ bounds a face, then $|E(P)|=1$, as otherwise the middle vertex of $P$ would contradict Lemma~\ref{lem:i012},
and $S(G)=\{k-2\}$, implying (\ref{cl-rep}).  Hence, we can assume $C_2$ does not bound a face, and thus $G_2$ is $C_2$-critical
by Lemma~\ref{lem:diskcritical}.  If $|E(P)|=2$, then $|C_2|=k$, and since $S(G)=S(G_2)$ and $|E(G_2)|<|E(G)|$, (\ref{cl-rep}) follows by the induction
hypothesis for $G_2$.  If $|E(P)|=1$, then $|C_2|=k-2$ and $S(G)=S(G_2)\in \Ss_{4,k-2}$, and thus (\ref{cl-rep}) holds.

Therefore, we can assume that $|C_1|\ge 5$ and by symmetry, $|C_2|\ge 5$.
If $|E(P)|=1$, then $|C_1|+|C_2|=k+2$, and thus $Z=\{|C_1|,|C_2|\}$ belongs to
$\Ss_{5,k}$ (the element corresponds to a ring of length $k$ with a chord).
Since $S(G_i)$ is either $\{|C_i|\}$ or belongs to $\Ss_{4,|C_i|}$ for $i\in\{1,2\}$, the multiset
$S(G)=S(G_1)\cup S(G_2)$ is a refinement of the element $Z$ of $\Ss_{5,k}$, implying (\ref{cl-rep}).
Finally, consider the case that $|E(P)|=2$.  As the middle vertex of $P$ has degree at least three, by symmetry we can assume that $C_2$
does not bound a face, and thus $G_2$ is $C_2$-critical.  In particular $|C_2|\ge 6$.
If $|C_2|=6$, then $S(G_2)=\emptyset$ by Theorem~\ref{thm:six} and $|C_1|=k-2$.  Consequently $S(G)=S(G_1)$
and either $S(G_1)=\{k-2\}$ or $S(G_1)$ belongs to $\Ss_{4,k-2}$, implying (\ref{cl-rep}).
If $|C_2|\ge 7$, then note that $|C_1|+(|C_2|-2)=k+2$, and thus $Z'=\{|C_1|,|C_2|-2\}$ belongs to $\Ss_{5,k}$.
Furthermore, $S(G_2)\in \Ss_{4,|C_2|}=\Ss_{4,(|C_2|-2)+2}$, and thus $\{C_1\}\cup S(G_2)$ is a one-step refinement of $Z'$.
Consequently, $S(G)$ is a refinement of the element $Z'$ of $\Ss_{5,k}$, implying (\ref{cl-rep}).

Therefore, we can assume $h_0$ is not ring-bound.  As $G$ is embedded in the disk, $G$ has no essential subgraph,
and thus Lemma~\ref{lemma-obstacle} applies; let $G'$ be the corresponding $R$-critical graph,
and let $J\subseteq G$ together with a system $\{S_f:f\in F(G')\}$ be the cover of $G$ by faces of $G'$.
By Lemma~\ref{lem:crit3conn}, each face of $G'$ is closed $2$-cell,
and thus each face of $G'$ has elasticity at most $2$.  

For each face $f\in F(G')$, the $G$-expansion of $S_f$ consists of $t=|S_f|$ graphs
$H_1$, \ldots, $H_t$ embedded in disks with rings $R_1$, \ldots, $R_t$ such that $|R_1|+\ldots+|R_t|=|f|+\el(f)$ and $|R_1|\ge \ldots\ge |R_t|$.
Let $m\in \{0,\ldots,t\}$ be the index such that $|R_i|>4$ for $1\le i\le m$ and $|R_i|=4$ for $m+1\le i\le t$.
Observe that if $t\ge 1+\el(f)/2$, then $\{|R_1|,\ldots, |R_m|\}$ is an element of $\Ss_{4,|f|}$
(corresponding to the graph obtained from a ring of length $|f|$ by adding $t-1$ chords to form faces
of lengths $|R_1|$, \ldots, $|R_{t-1}|$, and then adding further chords to split the last face
of length $|f|+2(t-1)-|R_1|-\ldots-|R_{t-1}|=|R_t|+2(t-1)-\el(f)$ into a face of length $|R_t|$ and $4$-faces),
and thus $S(H_1)\cup \ldots\cup S(H_t)$ is a
refinement of $\{|f|\}$.  If $t<1+\el(f)/2$, then $t=1$ and $\el(f)=2$ and Lemma~\ref{lemma-obstacle} asserts
that $S_f$ is non-trivial.  In that case, it follows that $H_1$ is $R_1$-critical and $|R_1|=|f|+2$, and
thus $S(H_1)$ belongs to $\Ss_{4,|f|+2}$; and again, $S(H_1)$ is a refinement of $\{|f|\}$.
Since this holds for every $f\in F(G')$ and $J$ with $\{S_f:f\in F(G')\}$ forms a cover of $G$ by faces of $G'$,
we conclude that $S(G)$ is a refinement of $S(G')$.

Since $|E(G')|<|E(G)|$, the induction hypothesis implies that either $S(G')=\{k-2\}$ or $S(G')$ is a refinement of an element of
$\Ss_{4,k-2}\cup \Ss_{5,k}$.  In the latter case, since $S(G)$ is a refinement of $S(G')$, it follows that
$S(G)$ is also a refinement of an element of $\Ss_{4,k-2}\cup \Ss_{5,k}$, and (\ref{cl-rep}) holds.
Hence, suppose that $S(G')=\{k-2\}$.  Let $f$ be the $(k-2)$-face of $G'$ and let $H_1$, \ldots, $H_t$ and $m$ be as in the
previous paragraph; we have $S(G)=S(H_1)\cup \ldots\cup S(H_m)$.  By the argument from the previous paragraph,
if $t\ge 1+\el(f)/2$, this implies $S(G)$ is a refinement of an element $\{|R_1|,\ldots, |R_m|\}$ of $\Ss_{4,|f|}=\Ss_{4,k-2}$.
Hence, suppose that $t=1$ and $\el(f)=2$, and thus $S_f$ is non-trivial.  Hence, $H_1$ is $R_1$-critical and $|R_1|=k$.
Note that $|E(H_1)|<|E(G)|$, since $R$ is a proper subgraph of $J$ by the definition of a cover by faces.
Since $S(G)=S(H_1)$, (\ref{cl-rep}) follows from the induction hypothesis applied to $H_1$.
\end{proof}

For a multiset $Z$, let $s(Z)=\sum_{z\in Z} s(z)$.

\begin{corollary}\label{cor-diskcase}
Let $G$ be a triangle-free graph embedded in the disk with one ring $R$ of length $l\ge 4$.  If $G$ is $R$-critical, then
$w(G,\{R\})\le s(l-2)$.
\end{corollary}
\begin{proof}
We proceed by induction on $l$.  The claim holds for $l\le 6$ by Theorem~\ref{thm:six}.
Suppose that $l\ge 7$.  By Lemma~\ref{lemma-diskcase}, we have that $S(G)$ is either $\{l-2\}$ or a refinement
of $Z\in\Ss_{4,l-2}\cup \Ss_{5,l}$.  In the former case, we have $w(G,\{R\})=s(l-2)$, and thus we can assume the latter.

If $Z\in\Ss_{4,l-2}$, then $s(Z)\le s(l-4)<s(l-2)$ by the induction hypothesis.  If $Z\in \Ss_{5,l}$, then $l\ge 8$ and
$s(Z)\le s(l-3)+s(5)<s(l-2)$ by Theorem~\ref{thm:diskgirth5}.
Since $S(G)$ is a refinement of $Z$, there exist multisets $Z=Z_0$, $Z_1$, \ldots, $Z_m=S(G)$ such that
$Z_{i+1}$ is a one-step refinement of $Z_i$ for $0\le i\le m-1$.  We claim that $s(Z_i)<s(l-2)$ for $i=1,\ldots,m$.
Indeed, assuming that the claim holds for $i-1$, the maximum of $Z_{i-1}$ is at most $l-3$,
and since $Z_i$ is a one-step refinement of $Z_{i-1}$, we have
$s(Z_{i})=s(Z_{i-1})-s(a)+s(A)$ for some $a\in Z_{i-1}$ and $A\in \Ss_{4,a}\cup \Ss_{4,a+2}$. Since $a+2<l$, by the induction hypothesis
we have $s(A)\le s((a+2)-2)=s(a)$, and thus $s(Z_i)\le s(Z_{i-1})<s(l-2)$.  We conclude that $S(G)=s(Z_m)<s(l-2)$.
\end{proof}

We will need the following properties of the function $\surf$, which were proved in~\cite{trfree3}.
\begin{lemma}\label{lemma-surfineq}
If $g$, $g'$, $t$, $t_0$, $t_1$, $t'_0$, $t'_1$ are non-negative integers, then
the following holds:
\begin{itemize}
\item[{\rm (a)}] Assume that if $g=0$ and $t\le 2$, then $t_0+t_1<t$.
If $t\ge 2$, $t'_0\le t_0$, $t'_1\le t_1$ and $t'_0+t'_1\ge t_0+t_1-2$, then
$\surf(g,t-1,t'_0, t'_1)\le \surf(g,t,t_0,t_1)-1$.
\item[{\rm (b)}] If $g' < g$ and either $g'>0$ or $t\ge 2$, then
$\surf(g',t,t_0, t_1)\le \surf(g,t,t_0, t_1)-120(g-g')+32$.
\item[{\rm (c)}] Let $g''$, $t'$, $t''$, $t''_0$ and $t''_1$
be nonnegative integers satisfying $g=g'+g''$, $t=t'+t''$, $t_0=t'_0+t''_0$, $t_1=t'_1+t''_1$,
either $g''>0$ or $t''\ge 1$,
and either $g'>0$ or $t'\ge 2$.  Then,
$\surf(g',t',t'_0, t'_1)+\surf(g'',t'',t''_0, t''_1)\le \surf(g,t,t_0, t_1)-\delta$,
where $\delta=16$ if $g''=0$ and $t''=1$ and $\delta=56$ otherwise.
\item[{\rm (d)}] If $g\ge 2$, then $\surf(g-2, t, t_0, t_1)\le \surf(g, t, t_0, t_1)-124$
\end{itemize}
\end{lemma}

Consider a graph $H$ embedded in a surface $\Pi$ with rings $\QQ$, and let $f$ be a face of $H$.
Let $a_0$ and $a_1$ be the number of weak and non-weak vertex-like rings, respectively, that form one of the facial walks of $f$ by
themselves.  Let $a$ be the number of facial walks of $f$.  We define $\surf(f)=\surf(g(\Pi_f),a,a_0,a_1)$.

Let $G_1$ be a graph embedded in $\Sigma_1$ with rings $\RR_1$ and $G_2$ a graph embedded in $\Sigma_2$ with rings $\RR_2$.
Let $m(G_i)$ denote the number of edges of $G_i$ that are not contained in the boundary of $\Sigma_i$.
Let us write $(G_1,\Sigma_1,\RR_1)\prec (G_2,\Sigma_2,\RR_2)$ to denote that 
the quadruple $(g(\Sigma_1),|\RR_1|,m(G_1),|E(G_1)|)$ is lexicographically smaller than $(g(\Sigma_2),|\RR_2|,m(G_2),|E(G_2)|)$.

We now have everything needed to prove Theorem~\ref{thm-maingen}.

\begin{proof}[Proof of Theorem~\ref{thm-maingen}]
Let $\eta_0$ be the constant from Theorem~\ref{thm:treti} and let $\eta=6\eta_0+44$.
We proceed by induction and assume that the claim holds for all graphs $G'$ embedded in surfaces $\Sigma'$ with rings $\RR'$
such that $(G',\Sigma',\RR')\prec (G,\Sigma,\RR)$.
Let $g=g(\Sigma)$, $t_0=t_0(\RR)$ and $t_1=t_1(\RR)$.
By Theorem~\ref{grotzsch}, Corollary~\ref{cor-diskcase} and Theorem~\ref{thm:gimtho}, the claim holds if $g+|\RR|\le 1$, hence assume that $g+|\RR|\ge 2$.
Similarly, if $g=0$ and $|\RR|=2$, then we can assume that $t_0+t_1\le 1$ by Corollary~\ref{cor-critshort}.
First, we aim to prove that we can assume that all faces in the embedding of $G$ are semi-closed $2$-cell.  For later use, we will consider a more general setting.

\claim{cl-non2nonomni}{Let $H$ be an internally triangle-free graph embedded in a surface $\Pi$ with rings $\QQ$ without non-contractible $4$-cycles, such that at least one face
of $H$ is not open $2$-cell and no face of $H$ is omnipresent.  If $H$ is $\QQ$-critical and $(H,\Pi,\QQ)\preceq (G,\Sigma,\RR)$, then
$$w(H,\QQ)\le \ell(\QQ) + \eta\Bigl(\surf(g(\Pi), |\QQ|, t_0(\QQ),t_1(\QQ)) - 7 - \sum_{h\in F(H)} \surf(h)\Bigr).$$}
\begin{proof}
We prove the claim by induction.
Consider for a moment an internally triangle-free graph $H'$ embedded in a surface $\Pi'$ with rings $\QQ'$ without non-contractible $4$-cycles, where
$(H',\Pi',\QQ')\prec(H,\Pi,\QQ)$, such that either $H'=\QQ'$ or $H'$ is $\QQ'$-critical.  We claim that
\begin{equation}\label{eq-ind}
w(H',\QQ')\le \ell(\QQ') + \eta\Bigl(\surf(g(\Pi'), |\QQ'|, t_0(\QQ'),t_1(\QQ')) - \sum_{h\in F(H')} \surf(h)\Bigr).
\end{equation}
The claim obviously holds if $H'=\QQ'$, since $w(H',\QQ')\le \ell(\QQ')$ in that case; hence, it suffices to consider the case that
$H'$ is $\QQ'$-critical.
If at least one face of $H'$ is not open $2$-cell and no face of $H'$ is omnipresent, then
this follows by an inductive application of (\ref{cl-non2nonomni}) (we could even strengthen the inequality by $7\eta$).
If all faces of $H'$ are open $2$-cell, then note that $\surf(h)=0$ for every $h\in H'$,
and since $(H',\Pi',\QQ')\prec(G,\Sigma,\RR)$, we can apply Theorem~\ref{thm-maingen} inductively to obtain (\ref{eq-ind}).
Finally, suppose that $H'$ has an omnipresent face $f$, let $\QQ'=\{Q_1,\ldots, Q_t\}$ and for $1\le i\le t$,
let $C_i$ be the cuff traced by $Q_i$, let $\Delta_i$ be a closed disk in $\Pi'+\widehat{C_i}$ such that $\widehat{C_i}\subset \Delta_i$ and the boundary of $\Delta_i$
is a subset of $f$, and let $f_i$ denote the boundary walk of $f$ contained in $\Delta_i$.
Since all components of $H'$ are planar and contain only one ring, Lemma~\ref{lem:crit3conn} implies that all
faces of $H'$ distinct from $f$ are closed $2$-cell.
Furthermore, by Theorem~\ref{thm:six}, each vertex-like ring forms component of the boundary of $f$ by itself, and hence
$\surf(f)=\surf(g(\Pi'),|\QQ'|,t_0(\QQ'),t_1(\QQ'))$.
If $Q_i$ is not a vertex-like ring, then by applying Corollary~\ref{cor-diskcase} to the subgraph $H'_i$ of $H'$ embedded in $\widehat{\Delta_i}\setminus \widehat{C_i}$, we conclude that the weight of $H'_i$ is at most
$s(|Q_i|)$ and that $|f_i|\le |Q_i|$.  Note that $s(|Q_i|)-s(|f_i|)\le |Q_i|-|f_i|$.  Therefore, we again obtain (\ref{eq-ind}):
\begin{align*}
w(H',\QQ')&\le |f|+\sum_{i=1}^t s(|Q_i|)-s(|f_i|)\le |f|+\sum_{i=1}^t |Q_i|-|f_i|\\
&=\ell(\QQ')\\
&=\ell(\QQ') + \eta\Bigl(\surf(g(\Pi'), |\QQ'|, t_0(\QQ'),t_1(\QQ')) - \sum_{h\in F(H')} \surf(h)\Bigr).
\end{align*}

Let us now return to the graph $H$.  Since $H$ is $\QQ$-critical, Theorem~\ref{grotzsch} implies that no component of $H$ is a planar graph without rings.
Let $f$ be a face of $H$ which is not open $2$-cell.  By the assumptions, $f$ is not omnipresent,
and thus $g(\Pi)>0$ or $|\QQ|>2$.  Let $c$ be a simple closed curve in $f$
infinitesimally close to a facial walk $W$ of $f$.
Cut $\Pi$ along $c$ and cap the resulting holes by disks ($c$ is always a $2$-sided curve).
Let $\Pi_1$ be the surface obtained this way that contains $W$, and if $c$ is separating, then let $\Pi_2$ be the other surface.
Since $f$ is not omnipresent, we can choose $W$ so that either $g(\Pi_1)>0$ or
$\Pi_1$ contains at least two rings of $\QQ$.  Let us discuss several cases:

\begin{itemize}
\item {\em $c$ is separating and $H$ is contained in $\Pi_1$.}  Let $H_1$ denote the embedding of $H$ in $\Pi_1$.
In this case $f$ has only one facial walk, and since $f$ is not open $2$-cell,
$\Pi_2$ is not the sphere.  It follows that $g(\Pi_1)=g(\Pi)-g(\Pi_2)<g(\Pi)$, and thus $(H_1,\Pi_1,\QQ)\prec (H,\Pi,\QQ)$.
Note that the weights of the faces of the embedding of $H$ in $\Pi$ and of $H_1$ in $\Pi_1$ are the same, with the exception of $f$
whose weight in $\Pi$ is $|f|$ and the weight of the corresponding face of $H_1$ in $\Pi_1$ is $s(|f|)\ge |f|-8$.
By (\ref{eq-ind}), we have
\begin{align*}
w(H,\QQ)&\le w(H_1,\QQ)+8\\
&\le\ell(\QQ) + 8+\eta\Bigl(\surf(g(\Pi_1), |\QQ|, t_0(\QQ),t_1(\QQ)) + \surf(f) - \sum_{h\in F(H)} \surf(h)\Bigr).
\end{align*}
Note that $\surf(f)=120g(\Pi_2)-72$.
By Lemma~\ref{lemma-surfineq}(b), we conclude that
$$w(H,\QQ)\le\ell(\QQ) +8 + \eta\Bigl(\surf(g(\Pi), |\QQ|, t_0(\QQ),t_1(\QQ)) -40 - \sum_{h\in F(H)} \surf(h)\Bigr).$$

\item {\em $c$ is separating and $\Pi_2$ contains a nonempty part $H_2$ of $H$.}  Let $H_1$ be the part of $H$ contained
in $\Pi_1$. Let $\QQ_i$ be the subset of $\QQ$ belonging to $\Pi_i$ and $f_i$ the face of $H_i$ corresponding to $f$, for $i\in \{1,2\}$.
Note that $f_1$ is an open disk, hence $\surf(f_1)=0$.  Using (\ref{eq-ind}), we get
\begin{align*}
w(H,\QQ)={}& w(f)-w(f_1)-w(f_2)+w(H_1,\QQ_1)+w(H_2,\QQ_2)\\
\le{}& w(f)-w(f_1)-w(f_2)+\ell(\QQ_1)+\ell(\QQ_2)+\\
&+\eta\sum_{i=1}^2\surf(g(\Pi_i), |\QQ_i|, t_0(\QQ_i),t_1(\QQ_i))+\\
&+\eta\Bigl(\surf(f)-\surf(f_2) - \sum_{h\in F(H)} \surf(h)\Bigr).
\end{align*}
Note that $w(f)-w(f_1)-w(f_2)\le 16$ and
$\ell(\QQ_1)+\ell(\QQ_2)=\ell(\QQ)$.  Also, $\surf(f)-\surf(f_2)\le 48$,
and when $g(\Pi_f)=0$ and $f$ has only two facial walks, then $\surf(f)-\surf(f_2)\le 6$.

By Lemma~\ref{lemma-surfineq}(c), we have 
$$\sum_{i=1}^2\surf(g(\Pi_i), |\QQ_i|, t_0(\QQ_i),t_1(\QQ_i)) \le\surf(g(\Pi), |\QQ|, t_0(\QQ),t_1(\QQ)) - \delta,$$
where $\delta=16$ if $g(\Pi_2)=0$ and $|\QQ_2|=1$ and $\delta=56$ otherwise.
Note that if $g(\Pi_2)=0$ and $|\QQ_2|=1$, then $g(\Pi_f)=0$ and $f$ has only two facial walks.
We conclude that $\surf(f)-\surf(f_2)-\delta \le -8$.  Therefore,
$$w(H,\QQ)\le \ell(\QQ)+16+\eta\Bigl(\surf(g(\Pi), |\QQ|, t_0(\QQ),t_1(\QQ))-8-\sum_{h\in F(H)} \surf(h)\Bigr).$$

\item {\em $c$ is not separating.} Let $f_1$ be the face of $H$ (in the embedding in $\Pi_1$) bounded by $W$ and $f_2$ the other face corresponding to
$f$.  Again, note that $\surf(f_1)=0$.  By (\ref{eq-ind}) applied to the embedding of $H$ in $\Pi_1$, we obtain the following for the weight of $H$ in $\Pi$:
\begin{align*}
w(H,\QQ)\le{}& w(f)-w(f_1)-w(f_2) + \ell(\QQ)+\\
&+ \eta\cdot\surf(g(\Pi_1), |\QQ|, t_0(\QQ),t_1(\QQ)) +\\
&+ \eta\Bigl(\surf(f) - \surf(f_2) - \sum_{h\in F(H)} \surf(h)\Bigr).
\end{align*}
Since $c$ is two-sided, $g(\Pi_1)=g(\Pi)-2$, and
$$\surf(g(\Pi_1), |\QQ|, t_0(\QQ),t_1(\QQ))=\surf(g(\Pi), |\QQ|, t_0(\QQ),t_1(\QQ))-124$$
by Lemma~\ref{lemma-surfineq}(d).
Since $\surf(f) - \surf(f_2)\le 48$ and $w(f)-w(f_1)-w(f_2)\le 16$, we have
$$w(H,\QQ)\le \ell(\QQ)+16+\eta\Bigl(\surf(g(\Pi), |\QQ|, t_0(\QQ),t_1(\QQ))-76-\sum_{h\in F(H)} \surf(h)\Bigr).$$
\end{itemize}
The results of all the subcases imply (\ref{cl-non2nonomni}).
\end{proof}

Next, we consider the case that one of the faces is omnipresent.

\claim{cl-omni}{Let $H$ be an internally triangle-free graph embedded in $\Sigma$ with rings $\RR$ 
and without non-contractible $4$-cycles, and let $f$ be an omnipresent face of $H$.
If $H$ is $\RR$-critical, then
$$w(H,\RR)\le \ell(\RR)-2=\ell(\RR) - 2 +  \eta\Bigl(\surf(g, |\RR|, t_0,t_1) - \sum_{h\in F(H)} \surf(h)\Bigr).$$}
\begin{proof}
Since $H$ is $\RR$-critical and $f$ is an omnipresent face, each component of $H$ is planar and contains
exactly one ring.  In particular, all faces of $H$ distinct from $f$ are closed $2$-cell.
For $R\in\RR$, let $H_R$ be the component of $H$ containing $R$.
Exactly one boundary walk $W$ of $f$ belongs to $H_R$.  Cutting along $W$ and capping the hole by a disk,
we obtain an embedding of $H_R$ in a disk with one ring $R$. Let $f_R$ be the face of this embedding bounded by $W$.
Note that either $H_R=R$ or $H_R$ is $R$-critical. If $R$ is a vertex-like ring, then we have $H_R=R$ by Theorem~\ref{thm:six}; hence,
every vertex-like ring in $\RR$ forms a facial walk of $f$, and thus $\surf(f)=\surf(g,|\RR|,t_0,t_1)$.
Consequently, $\surf(g, |\RR|, t_0,t_1)=\sum_{h\in F(H)} \surf(h)$, and it suffices to prove the first inequality of the claim.

Suppose that $H_R\neq R$ for a ring $R\in \RR$. By Corollary~\ref{cor-diskcase}, we have $w(H_R,\{R\})\le s(|R|-2)$.
Since $f_R$ is a face of $H_R$, it follows that $|f_R|\le |R|-2$.  Furthermore,
$w(H_R,\{R\})-w(f_R)\le s(|R|-2)-s(|f_R|)\le |R|-|f_R|-2$.  On the other hand, if $H_R=R$, then $w(H_R,\{R\})-w(f_R)=0=|R|-|f_R|$.
As $H$ is $\RR$-critical, there exists $R\in\RR$ such that $H_R\neq R$, and thus
\begin{align*}
w(H,\RR)&= w(f)+\sum_{R\in \RR} (w(H_R,\{R\})-w(f_R))\\
&\le|f|+\Bigl(\sum_{R\in\RR} |R|-|f_R|\Bigr)-2\\
&=\ell(\RR)-2.
\end{align*}
\end{proof}

Finally, consider the case that the embedding is open $2$-cell but not semi-closed $2$-cell.

\claim{cl-rep1}{Let $H$ be an internally triangle-free $\RR$-critical graph embedded in $\Sigma$ with rings $\RR$ (the same surface and rings as $G$)
so that all faces of $H$ are open $2$-cell and without non-contractible $4$-cycles.
If $|E(H)|\le |E(G)|$ and a face $f$ of $H$ is not semi-closed $2$-cell,
then $$w(H,\RR)\le \ell(\RR) + \eta\Bigl(\surf(g,|\RR|, t_0, t_1)-1/2\Bigr).$$}
\begin{proof}
Since $f$ is not semi-closed $2$-cell, there exists a vertex $v$ appearing at least twice in the facial walk of $f$,
and a simple closed curve $c$ going through the interior of $f$ and joining two of the appearances of $v$;
furthermore, if $v$ is the main vertex of a vertex-like ring, it is not the case that $c$ separates $R$ from the rest of the graph $H$.
Cut the surface along $c$ and patch the resulting hole(s) by disk(s).
Let $v_1$ and $v_2$ be the two vertices to that $v$ is split.
For $i=1,2$, if $v_i$ is not incident with a cuff in the resulting surface, drill a new hole next to it in the incident patch and add a triangle $T_i$ tracing its boundary,
with vertex set consisting of $v_i$ and two new vertices (note that if $v$ is an internal vertex, a hole will be drilled at both
$v_1$ and $v_2$, while if $v$ is a ring vertex, a hole will be drilled at exactly one of $v_1$ and $v_2$).

If the curve $c$ is separating, then let $H_1$ and $H_2$ be the resulting graphs embedded in the two
surfaces $\Sigma_1$ and $\Sigma_2$ obtained by this construction; if $c$ is not separating, then let $H_1$
be the resulting graph embedded in a surface $\Sigma_1$.  We choose the labels so that $v_1\in V(H_1)$.
If $c$ is two-sided, then let $f_1$ and $f_2$ be the faces to that $f$ is split by $c$, where $f_1$ is a face of $H_1$.
If $c$ is one-sided, then let $f_1$ be the face in $\Sigma_1$ corresponding to $f$.
Note that $|f_1|+|f_2|>|f|$ in the former case, and thus
$w(f)-w(f_1)-w(f_2)<16$.  In the latter case, we have $w(f)=w(f_1)$.

If $c$ is separating, then for $i\in\{1,2\}$, let $\RR_i$ consist of the rings of $\RR\setminus \{v\}$ contained in $\Sigma_i$, and
if none of these rings contains $v_i$ (so that $T_i$ exists), then also of the vertex-like ring $T_i$.
Here, we designate $T_i$ as weak if $v$ is an internal vertex, $\Sigma_{3-i}$ is a cylinder and the ring of $H_{3-i}$ distinct from $T_{3-i}$ is a vertex-like ring.
If $c$ is not separating, then let $\RR_1$ consist of the rings of $\RR$, together with those of $T_1$ and $T_2$ that exist.
In this case, we treat $T_1$ and $T_2$ as non-weak vertex-like rings.

Suppose first that $c$ is not separating.  By Lemma~\ref{lemma-crcon}, $H_1$ is $\RR_1$-critical.  Note that $H_1$ has at most two more (non-weak vertex-like) rings than $H$ and
$g(\Sigma_1)\in \{g-1,g-2\}$ (depending on whether $c$ is one-sided or not), and that $H_1$ has at least two rings.
If $H_1$ has only one more ring than $H$, then
\begin{align*}
\surf(g(\Sigma_1),|\RR_1|,t_0(\RR_1),t_1(\RR_1))&\le \surf(g-1,|\RR|+1,t_0,t_1+1)\\
&\le\gen(g-1,|\RR|+1,t_0,t_1+1)+32\\
&=\gen(g,|\RR|,t_0,t_1)-44\\
&=\surf(g,|\RR|,t_0,t_1)-44.
\end{align*}
Let us now consider the case that $H_1$ has two more rings than $H$ (i.e., that $v$ is an internal vertex).  If $g(\Sigma_1)=0$ and $|\RR_1|=2$,
then note that both rings of $H_1$ are vertex-like rings.  Corollary~\ref{cor-critshort} implies that
$H_1$ has only one edge; but the corresponding edge in $H$ would form a loop, which is a contradiction.
Consequently, we have $g(\Sigma_1)>1$ or $|\RR_1|\ge 3$, and
\begin{align*}
\surf(g(\Sigma_1),|\RR_1|,t_0(\RR_1),t_1(\RR_1))&\le \surf(g-1,|\RR|+2,t_0,t_1+2)\\
&=\surf(g,|\RR|,t_0,t_1)-32.
\end{align*}
By the induction hypothesis, we can apply Theorem~\ref{thm-maingen} to $H_1$, concluding that
$w(H,\RR)\le \ell(\RR)+22+\eta\Bigl(\surf(g,|\RR|,t_0,t_1) - 32\Bigr)$, and the claim follows.

Next, we consider the case that $c$ is separating.  Let us remark that $H_i$ is $\RR_i$-critical for $i\in\{1,2\}$.
This follows from Lemma~\ref{lemma-crcon} unless $T_i$ is a weak vertex-like ring.  If $T_i$ is a weak vertex-like ring, then $\Sigma_{3-i}$ is a cylinder
and the ring $R_{3-i}$ of $H_{3-i}$ distinct from $T_{3-i}$ is a vertex-like ring.  By Corollary~\ref{cor-critshort},
$H_{3-i}$ is a single edge and $R_{3-i}$ is not weak.  Consider any edge $e$ of $H_i$ not belonging to the rings.  Since $H$ is $\RR$-critical,
there exists a precoloring $\psi$ of $\RR$ that extends to $H-e$, but not to $H$.  Let $\psi_1$ be obtained from $\psi$ by giving a non-main vertex
of $T_i$ the color of the main vertex of $R_{3-i}$.  Then clearly $\psi_1$ extends to $H_i-e$, but not to $H$.  This implies that $H_i$ is indeed $\RR_i$-critical.

Thus, we can apply Theorem~\ref{thm-maingen} inductively for $H_1$ and $H_2$, and we have
\begin{align*}
w(H,\RR)&=w(H_1,\RR_1)+w(H_2,\RR_2) + w(f)-w(f_1)-w(f_2)\\
&\le \ell(\RR)+22 + \eta\sum_{i=1}^2\surf(g(\Sigma_i),|\RR_i|,t_0(\RR_i),t_1(\RR_i))
\end{align*}
Therefore, since $\eta\ge 44$, to prove (\ref{cl-rep1}) it suffices to prove that
\begin{equation}\label{eq-loccut}
\sum_{i=1}^2\surf(g(\Sigma_i),|\RR_i|,t_0(\RR_i),t_1(\RR_i))\le \surf(g,|\RR|,t_0,t_1)-1.
\end{equation}

Let us first consider the case that say $g(\Sigma_1)=0$ and $|\RR_1|=1$; i.e., $c$ passes through a vertex incident with a ring $R$ of $H$
and surrounds this ring.  By Lemma~\ref{lem:crit3conn}, the graph $H_1$ is $3$-connected, and in particular by the choice of $c$,
we conclude the ring $R$ is not vertex-like.  Hence, $\surf(g(\Sigma_1), |\RR_1|,t_0(\Sigma_1), t_1(\Sigma_1))=0$,
$g(\Sigma_2)=g(\Sigma)$, $|\RR_2|=|\RR|$, $t_0(\RR_2)=t_0(\RR)$ and $t_1(\RR_2)=t_1(\RR)+1$.
This implies $$\sum_{i=1}^2\surf(g(\Sigma_i),|\RR_i|,t_0(\RR_i),t_1(\RR_i))=\surf(g,|\RR|,t_0,t_1)-\delta,$$
where $\delta=8$ if $g(\Sigma_2)=0$, $|\RR_2|=2$, and both rings of $\RR_2$ are vertex-like,
and $\delta=4$ otherwise.

Therefore, we can assume that for $i\in\{1,2\}$, if $g(\Sigma_i)=0$, then $|\RR_i|\ge 2$.
If $|\RR_1|+|\RR_2|=|\RR|+1$, we have
\begin{align*}
\sum_{i=1}^2\surf(g(\Sigma_i),|\RR_i|,t_0(\RR_i),t_1(\RR_i))&\le \sum_{i=1}^2(\gen(g(\Sigma_i),|\RR_i|,t_0(\RR_i),t_1(\RR_i))+32)\\
&=\gen(g,|\RR|,t_0,t_1)-12\\
&=\surf(g,|\RR|,t_0,t_1)-12.
\end{align*}

Therefore, we can assume that $|\RR_1|+|\RR_2|=|\RR|+2$, i.e., $v$ is an internal vertex.
If for both $i\in\{1,2\}$, we have $g(\Sigma_i)>0$ or $|\RR_i|>2$, then
$g(\Sigma_1)+g(\Sigma_2)=g(\Sigma)$, $t_0(\RR_1)+t_0(\RR_2)=t_0(\RR)$, $t_1(\RR_1)+t_1(\RR_2)=t_1(\RR)+2$, and
\begin{align*}
\sum_{i=1}^2\surf(g(\Sigma_i),|\RR_i|,t_0(\RR_i),t_1(\RR_i))&= \sum_{i=1}^2\gen(g(\Sigma_i),|\RR_i|,t_0(\RR_i),t_1(\RR_i))\\
&=\gen(g,|\RR|,t_0,t_1)-32\\
&=\surf(g,|\RR|,t_0,t_1)-32.
\end{align*}
and the claim follows.

Hence, we can assume that say $g(\Sigma_1)=0$ and $|\RR_1|=2$.  Then, $\RR_1=\{T_1,R_1\}$ for some ring $R_1$,
$g(\Sigma_2)=g$ and $|\RR_2|=|\RR|$.  Since $H_1$ is $\RR_1$-critical, Corollary~\ref{cor-critshort}
implies that $R_1$ is not a weak vertex-like ring.  If $R_1$ is a vertex-like ring, then $T_2$ is a weak
vertex-like ring of $\RR_2$ which replaces the non-weak vertex-like ring $R_1$.
Therefore,
$\surf(g(\RR_2),|\RR_2|,t_0(\RR_2),t_1(\RR_2))=\surf(g,|\RR|,t_0,t_1)-1$.
Furthermore, $\surf(g(\RR_1),|\RR_1|,t_0(\RR_1),t_1(\RR_1))=\surf(0,2,0,2)=0$,
and the claim follows.

Finally, consider the case that $R_1$ is a non-vertex-like ring.  By symmetry, we can assume that
if $g(\Sigma_2)=0$ and $|\RR_2|=2$, then $\RR_2$ also contains a non-vertex-like ring.
Since $\RR_2$ is obtained from $\RR$ by replacing $R_1$ by a non-weak vertex-like ring $T_2$,
we have $\surf(g(\RR_2),|\RR_2|,t_0(\RR_2),t_1(\RR_2))=\surf(g,|\RR|,t_0,t_1)-4$.
Furthermore, $\surf(g(\RR_1),|\RR_1|,t_0(\RR_1),t_1(\RR_1))=\surf(0,2,0,1)=2$.
Consequently,
$$\sum_{i=1}^2\surf(g(\Sigma_i),|\RR_i|,t_0(\RR_i),t_1(\RR_i))\le \surf(g,|\RR|,t_0,t_1)-2.$$

Therefore, the inequality (\ref{eq-loccut}) holds.  This proves (\ref{cl-rep1}).
\end{proof}

By (\ref{cl-non2nonomni}), (\ref{cl-omni}) and (\ref{cl-rep1}), we either have $w(G,{\RR})<\eta\cdot\surf(g(\Sigma),|\RR|, t_0(\RR), t_1(\RR)) + \ell({\RR})$
or every face of $G$ is semi-closed $2$-cell.  From now on, assume that the latter holds.

Suppose that there exists a path $P\subset G$ of length at most $12$ with ends in distinct rings $R_1,R_2\in {\RR}$.
By choosing the shortest such path, we can assume that $P$ intersects no other rings.
Let $J=P\cup \bigcup_{R\in {\RR}} R$ and let $S=\{f\}$, where $f$ is the unique face of $J$.
Let $\{G'\}$ be the $G$-expansion of $S$, let $\Sigma'$ be the surface in that $G'$ is embedded and let ${\RR}'$ be the natural rings of $G'$.
Note that $g(\Sigma')=g$, $|\RR'|=|\RR|-1$, $\ell(\RR')\le \ell(\RR)+24$ and $t_0(\RR')+t_1(\RR')\ge t_0+t_1-2$.
Since $(G',\Sigma',\RR')\prec (G,\Sigma,\RR)$, by induction and by Lemma~\ref{lemma-surfineq}(a) we have
$w(G,{\RR})=w(G',{\RR'})\le \eta\cdot\surf(g,|\RR|-1,t_0(\RR'),t_1(\RR'))+\ell(\RR)+24<\eta\cdot\surf(g,|\RR|,t_0,t_1)+\ell(\RR)$.
Therefore, we can assume that no such path exists in~$G$.

Let $R\in\RR$ be a ring of $G$.  Suppose $G$ contains a subgraph $P$ with at most $12$ edges such that $P$ is
\begin{itemize}
\item a path joining two distinct vertices $u$ and $v$ of $R$ and otherwise disjoint from $R$ such that none of the cycles in $R\cup P$
is contractible, or
\item a non-contractible cycle not surrounding $R$, containing a vertex $u=v$ of $R$ and otherwise disjoint from $R$, or
\item a non-contractible cycle not surrounding $R$ and disjoint from $R$, together with a path joining one of its vertices
to a vertex $u=v$ of $R$.
\end{itemize}
Note that the existence of $P$ implies that either $g>0$ or $|\RR|\ge 3$, and thus $\surf(g,|\RR|,t_0(\RR),t_1(\RR))=\gen(g,|\RR|,t_0(\RR),t_1(\RR))$.
Since the distance between any two rings in $G$ is at least $13$, all vertices of $V(P)\setminus \{u,v\}$ are internal.
Let $J$ be the subgraph of $G$ consisting of $P$ and of the union of the rings, and let $S$ be the set of all faces of $J$.
Let $\{G_1,\ldots, G_k\}$ be the $G$-expansion of $S$,
and for $1\le i\le k$, let $\Sigma_i$ be the surface in that $G_i$ is embedded and let $\RR_i$ be the natural rings of $G_i$.  Note that
$\sum_{i=1}^k t_0(\RR_i)=t_0-b_0$ and $\sum_{i=1}^k t_1(\RR_i)=t_1-b_1$, where
$b_0=1$ if $R$ is a weak vertex-like ring and $b_0=0$ otherwise,
and $b_1=1$ if $R$ is a non-weak vertex-like ring and $b_1=0$ otherwise.
Let $r=\left(\sum_{i=1}^k |\RR_i|\right)-|\RR|$ and observe that
either $r=0$ and $k=1$, or $r=1$ and $1\le k\le 2$ (depending on whether the curves in $\widehat{\Sigma}$ corresponding to cycles in $R\cup P$
distinct from $R$ are one-sided, two-sided and non-separating or two-sided and separating). Furthermore, as seen from this case analysis,
$\sum_{i=1}^k g(\Sigma_i)=g+2k-r-3$.

We claim that $(G_i,\Sigma_i,\RR_i)\prec(G,\Sigma,\RR)$ for $1\le i\le k$.  This is clearly the case, unless $g(\Sigma_i)=g$.
Then, we have $k=2$, $r=1$ and $g(\Sigma_{3-i})=0$.  Since the boundaries of faces of $R\cup P$ are non-contractible, $\Sigma_{3-i}$ is not a disk,
hence $|\RR_{3-i}|\ge 2$ and $|\RR_i|<|\RR|$, again implying $(G_i,\Sigma_i,\RR_i)\prec(G,\Sigma,\RR)$.

By induction, we have
$w(G_i,\RR_i)\le \ell(\RR_i)+\eta\cdot\surf(g(\Sigma_i),|\RR_i|,t_0(\RR_i),t_1(\RR_i))$, for $1\le i\le k$.
Since every face of $G$ corresponds to a face of $G_i$ for some $i\in\{1,\ldots, k\}$
and $\sum_{i=1}^k\ell(\RR_i)\le \ell(\RR)+24$, we conclude that
$$w(G,\RR)\le \ell(\RR)+24+\eta\sum_{i=1}^k\surf(g(\Sigma_i),|\RR_i|,t_0(\RR_i),t_1(\RR_i)).$$
Note that for $1\le i\le k$, we have that $\Sigma_i$ is not a disk and $\RR_i$ contains at least one non-vertex-like ring,
and thus $\surf(g(\Sigma_i),|\RR_i|,t_0(\RR_i),t_1(\RR_i))\le \gen(g(\Sigma_i),|\RR_i|,t_0(\RR_i),t_1(\RR_i))+30$.
Therefore,
\begin{align*}
\lefteqn{\sum_{i=1}^k\surf(g(\Sigma_i),|\RR_i|,t_0(\RR_i),t_1(\RR_i))}\\
&\le\sum_{i=1}^k(\gen(g(\Sigma_i),|\RR_i|,t_0(\RR_i),t_1(\RR_i))+30)\\
&=\gen(g,|\RR|,t_0,t_1)+120(2k-r-3)+48r-120(k-1)+30k+4b_1+5b_0\\
&=\gen(g,|\RR|,t_0,t_1)+150k-72r-240+4b_1+5b_0\\
&\le\gen(g,|\RR|,t_0,t_1)-7=\surf(g,|\RR|,t_0,t_1)-7.
\end{align*}
The inequality of Theorem~\ref{thm-maingen} follows.  Therefore, we can assume that $G$ contains no such subgraph $P$.

Suppose that $G$ contains a connected essential subgraph $H$ with at most $12$ edges
(implying that $g>0$ or $|\RR|\ge 3$, and that $\surf(g,|\RR|,t_0(\RR),t_1(\RR))=\gen(g,|\RR|,t_0(\RR),t_1(\RR))$).
We can assume that $H$ is minimal,
i.e., no proper connected subgraph of $H$ is essential.
Given the already excluded cases, all vertices of $H$ are internal and $H$ satisfies one of the following conditions (see~\cite{rs7}):
\begin{itemize}
\item[(a)] $H$ is a non-contractible cycle not surrounding any of the rings, or
\item[(b)] $H$ is the union of two cycles $C_1$ and $C_2$ intersecting in exactly one vertex, and $C_1$ and $C_2$ surround distinct rings, or
\item[(c)] $H$ is the union of two vertex-disjoint cycles $C_1$ and $C_2$ and of a path between them, and $C_1$ and $C_2$ surround distinct rings, or
\item[(d)] $H$ is the theta-graph and each of the three cycles in $G$ surrounds a different ring (in this case, $\Sigma$ is the sphere with three holes).
\end{itemize}
Let $J$ be the subgraph of $G$ consisting of $H$ and of the union of the rings, and let $S$ be the set of all faces of $J$.
Let $\{G_1,\ldots, G_k\}$ be the $G$-expansion of $S$,
and for $1\le i\le k$, let $\Sigma_i$ be the surface in that $G_i$ is embedded and let $\RR_i$ be the natural rings of $G_i$.  Let $r=\left(\sum_{i=1}^k |\RR_i|\right)-|\RR|$.
Note that $\sum_{i=1}^k t_0(\RR_i)=t_0$, $\sum_{i=1}^k t_1(\RR_i)=t_1$ and $\sum_{i=1}^k\ell(\RR_i)\le \ell(\RR)+24$.

Let us first consider the case (a).  Observe that we have either $r=1$ and $k=1$, or $r=2$ and $1\le k\le 2$, and that $\sum_{i=1}^k g(\Sigma_i)=g-r+2k-2$.
If $g(\Sigma_1)=g$, then $k=2$ and $g(\Sigma_2)=0$; furthermore, $\Sigma_2$ has at least three cuffs, since $H$ does not surround a ring.
Thus, if $g(\Sigma_1)=g$, then $r=2$ and 
$|\RR_1|=|\RR|+r-|\RR_2|<|\RR|$.  The same argument
can be applied to $\Sigma_2$ if $k=2$, hence $(G_i,\Sigma_i,\RR_i)\prec(G,\Sigma,\RR)$ for $1\le i\le k$.  

By induction, we conclude that
\begin{equation}\label{eq-indh}
w(G,\RR)\le \ell(\RR)+24+\eta\sum_{i=1}^k\surf(g(\Sigma_i),|\RR_i|,t_0(\RR_i),t_1(\RR_i)).
\end{equation}
For $1\le i\le k$, let $\delta_i=72$ if  $g(\Sigma_i)=0$ and $|\RR_i|=1$, let $\delta_i=30$ if $g(\Sigma_i)=0$ and $|\RR_i|=2$, and let $\delta_i=0$ otherwise,
and note that since $\RR_i$ contains a non-vertex-like ring, we have
$$\surf(g(\Sigma_i),|\RR_i|,t_0(\RR_i),t_1(\RR_i))=\gen(g(\Sigma_i),|\RR_i|,t_0(\RR_i),t_1(\RR_i))+\delta_i.$$

If $k=2$, then recall that since $H$ does not surround a ring, we have either $g(\Sigma_i)>0$ or $|\RR_i|\ge 3$ for
$i\in \{1,2\}$; hence, $\delta_1+\delta_2=0$.
If $k=1$, then recall that we can assume that $G$ is not embedded in the projective plane with no rings by Theorem~\ref{thm:gimtho};
hence, either $g(\Sigma_1)>0$, or $|\RR_1|\ge 2$.  Consequently, we have $\delta_1\le 30$.

Combining the inequalities, we obtain $\sum_{i=1}^k \delta_i\le 60-30k$, and
\begin{align*}
\lefteqn{\sum_{i=1}^k\surf(g(\Sigma_i),|\RR_i|,t_0(\RR_i),t_1(\RR_i))}\\
&=\sum_{i=1}^k\gen(g(\Sigma_i),|\RR_i|,t_0(\RR_i),t_1(\RR_i))+\delta_i\\
&\le\gen(g,|\RR|,t_0,t_1)+120(2k-r-2)+48r-120(k-1)+\sum_{i=1}^k\delta_i\\
&\le\gen(g,|\RR|,t_0,t_1)+90k-72r-60\\
&\le\surf(g,|\RR|,t_0,t_1)-24.
\end{align*}
Together with (\ref{eq-indh}), this implies the inequality of Theorem~\ref{thm-maingen}.

If $H$ satisfies (b), (c) or (d), then $k=3$, $r=3$ and say $\Sigma_1$ and $\Sigma_2$ are cylinders and $\Sigma_3$ is obtained from $\Sigma$ by replacing two cuffs by one.
Therefore, we can again apply the induction hypothesis to obtain (\ref{eq-indh}).  Furthermore, denoting $b_{i,j}=t_j(\RR_i)$
(so $0\le b_{i,0}+b_{i,1}\le 1$) for $i\in\{1,2\}$ and $j\in\{0,1\}$, 
\begin{align*}
\lefteqn{\sum_{i=1}^k\surf(g(\Sigma_i),|\RR_i|,t_0(\RR_i),t_1(\RR_i))}\\
&=\surf(g,|\RR|-1,t_0-b_{1,0}-b_{2,0},t_1-b_{1,1}-b_{2,1})+\sum_{i=1}^2\surf(0,2,b_{i,0},b_{i,1})\\
&\le\surf(g,|\RR|,t_0,t_1)-6.
\end{align*}
This again implies the inequality of Theorem~\ref{thm-maingen}.
Therefore, assume that 
\claim{cl-essential}{\qquad every connected essential subgraph of $G$ has at least $13$ edges.}

Suppose now $G$ contains a path $P$ of length one or two intersecting $\bigcup\RR$ exactly in its endpoints $u$ and $v$.
Both ends of $P$ belong to the same ring $R$; let $P$, $P_1$ and $P_2$ be the paths in $R\cup P$ joining $u$ and $v$.
Recall that we can assume that $P\cup P_2$ is a contractible cycle.
Let $J$ be the subgraph of $G$ consisting of $P$ and of the union of the rings, and let $S$ be the set of faces of $J$.
Let $\{G_1,G_2\}$ be the $G$-expansion of $S$,
and for $1\le i\le 2$, let $\Sigma_i$ be the surface in that $G_i$ is embedded and let $\RR_i$ be the natural rings of $G_i$.  Note that say $\Sigma_2$ is a disk and its ring has length $|P|+|P_2|$,
and $\Sigma_1$ is homeomorphic to $\Sigma$.  Therefore, by induction and Corollary~\ref{cor-diskcase}, we have
$w(G,\RR)=w(G_1,\RR_1)+w(G_2,\RR_2)\le \ell(\RR_1)+\surf(g,|\RR|,t_0,t_1)+s(\ell(\RR_2))$.
Furthermore, $\ell(\RR_1)+s(\ell(\RR_2))\le \ell(\RR_1)+\ell(\RR_2)-4=\ell(R)+2|P|-4\le \ell(R)$, and the claim of the theorem follows.
Therefore, we can assume that
\claim{cl-nopath}{$G$ contains no path of length one or two intersecting $\RR$ exactly in its endpoints.}

Note that every $4$-cycle in $G$ is contractible, and thus it bounds a face by Lemma~\ref{lem:diskcritical} and Theorem~\ref{thm:six}.
Let us consider the case that every $4$-face in $G$ is ring-bound.  Since the distance between rings is at least $13$,
by (\ref{cl-nopath}) this implies that every $4$-face is incident with the main vertex of a vertex-like ring.  
For each vertex-like ring $R\in\RR$ whose main vertex $r$ is incident with a $4$-face, choose one such $4$-face bounded by a cycle
$rw_2w_3w_4$ and let $J_R$ be the edge $rw_2$ together with the $(\le\!6)$-cycle surrounding $R$
and containing the path $w_2w_3w_4$.  For any other ring $R\in\RR$, let $J_R=\emptyset$.
Since the distance between rings is at least $13$, the subgraphs $J_R$ for $R\in\RR$ are pairwise vertex-disjoint.
Let $J=\bigcup_{R\in\RR} R\cup J_R$
and let $S$ be the set of all faces of $J$.  Let $\{G_1,\ldots, G_k\}$ be the $G$-expansion of $S$,
and for $1\le i\le k$, let $\Sigma_i$ be the surface in that $G_i$ is embedded and let $\RR_i$ be the natural rings of $G_i$, labelled so that $\Sigma_1$ is homeomorphic to $\Sigma$
and $\Sigma_2$, \ldots, $\Sigma_k$ are disks bounded by rings of length at most $11$, corresponding to the vertex-like rings
of $G$ incident with $4$-faces.  Note that $G_1$ has internal girth at least five, and thus
$w(G_1,\RR_1)\le \ell(R_1)+\eta_0\cdot\surf(g,|\RR|,0,0)$ by Theorem~\ref{thm:treti}.  Furthermore, $w(G_i,\RR_i)\le s(9)$ for $2\le i\le k$
by Corollary~\ref{cor-diskcase}.  Therefore,
\begin{align*}
w(G,\RR)&=\sum_{i=1}^k w(G_i,\RR_i)\\
&\le \ell(\RR_1)+\eta_0\cdot\surf(g,|\RR|,0,0)+|\RR|s(9)\\
&\le \ell(\RR) + \eta_0\cdot\surf(g,|\RR|,0,0)+|\RR|(6+s(9))\\
&\le \ell(\RR) + \eta\cdot\surf(g,|\RR|,t_0,t_1).
\end{align*}

Finally, suppose that $G$ contains a $4$-face which is not ring-bound.
Let $G'$ be the $\RR$-critical graph embedded in $\Sigma$ such that $|E(G')|<|E(G)|$,
and let $J\subseteq G$ and $\{S_f:f\in F(G')\}$ be the cover of $G$ by faces of $G'$, obtained by Lemma~\ref{lemma-obstacle}.
If $G'$ does not contain any non-contractible $4$-cycle, then let $G_0=G'$ and $\RR_0=\RR$.

Otherwise, consider the unique non-contractible $4$-cycle $C=w_1w_2w_3w_4$ in $G'$, which is flippable to a $4$-face, and
let $R\in \RR$ be the ring surrounded by $C$.  Let $\Delta$ be the closed disk bounded by $C$ in $\Sigma+\widehat{C}_R$,
where $C_R$ is the cuff corresponding to $R$.  By symmetry, we can assume that a face $f_1$ whose boundary contains
the path $w_1w_2w_3$ is a subset of $\Delta$, and the face $f_2$ whose boundary contains
the path $w_1w_4w_3$ is not a subset of $\Delta$. See Figure~\ref{fig-flip}, which also illustrates the transformation
we describe in the following paragraph.

Note that the subgraph $H$ of $G'$ drawn in $\Delta$ is connected (since $G'$ is critical, Theorem~\ref{grotzsch} implies that
any component of $H$ contains $R$ or $C$, and if $C$ were in a different component from $R$, then by Theorem~\ref{thm:six}
the cycle $C$ would form a component of $H$ by itself, and due to the face $f_2$, $w_4$ would be an internal vertex of $G'$
of degree two, contradicting Lemma~\ref{lem:i012}).
Hence, the face $f_1$ is open $2$-cell.
Let $a$ denote the walk such that the boundary of $f_1$ is the concatenation of $a$ and $w_1w_2w_3$.
Note that $w_1$ and $w_3$ form a $2$-cut in $G'$.  Let $G_0$ be the embedding of $G'$ in $\Sigma$ obtained by mirroring $H-w_2$
(i.e., the cyclic orders of neighbors of vertices in $H-\{w_1,w_2,w_3\}$ are reversed, and the orders of the neighbors of $w_1$ and $w_3$
in $H-w_2$ are reversed).  Let $\RR_0$ be the rings of $G_0$ obtained from $\RR$ by mirroring $R$.
Note that $G_0$ contains no non-contractible $4$-cycle, as $C$ is the only $4$-cycle in $G'$ containing $w_1$ and $w_3$
and $C$ becomes a $4$-face in $G_0$.
Let $f$ be the face of $G_0$ whose boundary contains the walk $a$.  Note that $\Sigma_f=\Sigma_{f_2}$.
We have $w(G',\RR)=w(G_0,\RR)-w(f)+w(f_1)+w(f_2)$.  If $f_2$ is not open $2$-cell, then
$w(f)-w(f_1)-w(f_2)=|f|-s(|f_1|)-|f_2|=|f_1|-4-s(|f_1|)\ge 0$.  If $f_2$ is open $2$-cell,
then $w(f)-w(f_1)-w(f_2)=s(|f|)-s(|f_1|)-s(|f_2|)=s(|f_1|+|f_2|-4)-s(|f_1|)-s(|f_2|)\ge 0$.
Therefore,
\begin{equation}\label{ineq-flip}
w(G',\RR)\le w(G_0,\RR_0).
\end{equation}
Note that (\ref{ineq-flip}) holds trivially when $G_0=G'$.

For $f\in F(G')$, let $\{G^f_1, \ldots, G^f_{k_f}\}$
be the $G$-expansion of $S_f$ and for $1\le i\le k_f$, let $\Sigma^f_i$ be the surface in that $G^f_i$ is embedded and let $\RR^f_i$ denote the natural rings of $G^f_i$.
Since $J$ and $\{S_f:f\in F(G')\}$ is a cover of $G$ by faces of $G'$, we have 
\begin{equation}\label{eq-cover}
w(G,\RR)=\sum_{f\in F(G')} \sum_{i=1}^{k_f} w(G^f_i,\RR^f_i).
\end{equation}

Consider a face $f\in F(G')$.  We have $g(\Sigma_f)\le g$.  If $g(\Sigma_f)=g$, then every component of $G'$ is planar,
and since $G'$ is $\RR$-critical and triangle-free, Theorem~\ref{grotzsch} implies that each component of $G'$ contains
at least one ring of $\RR$; consequently, $f$ has at most
$|\RR|$ facial walks and $\Sigma_f$ has at most $|\RR|$ cuffs.  Since the surfaces embedding the components of the $G$-expansion of $S_f$ are
fragments of $\Sigma_f$, we have $(G^f_i,\Sigma^f_i,\RR^f_i)\prec (G,\Sigma,\RR)$ for $1\le i\le k_f$:
otherwise, we would have $m(G^f_i)=m(G)$, hence by the definition of $G$-expansion, the boundary of $S_f$ would have to be equal
to the union of rings in $\RR$, contrary to the fact that $J\neq \bigcup \RR$ by the definition of a cover of $G$ by faces of $G'$.

Therefore, we can apply Theorem~\ref{thm-maingen}
inductively for $G^f_i$ and we get 
\begin{equation}\label{ineq-fromin}
w(G^f_i,\RR^f_i)\le \ell(\RR^f_i)+\eta\cdot\surf(g(\Sigma^f_i),|\RR^f_i|,t_0(\RR^f_i), t_1(\RR^f_i)).
\end{equation}
Observe that since $\{\Sigma^f_1,\ldots,\Sigma^f_{k_f}\}$ are fragments of $\Sigma_f$, we have
$$\sum_{i=1}^{k_f} \surf(g(\Sigma^f_i),|\RR^f_i|,t_0(\RR^f_i),t_1(\RR^f_i))\le \surf(f),$$
and summing (\ref{ineq-fromin}) for $i=1,\ldots, k_f$, we obtain
\begin{equation}\label{eq-partsfprel}
\sum_{i=1}^{k_f} w(G^f_i,\RR^f_i)\le |f|+\el(f)+\eta\cdot\surf(f).
\end{equation}
In case that $f$ is semi-closed $2$-cell, all fragments of $f$ are disks and we can use Corollary~\ref{cor-diskcase}
instead of Theorem~\ref{thm-maingen} and strengthen the inequality (\ref{ineq-fromin}) to
$$w(G^f_i,\RR^f_i)\le s(\ell(\RR^f_i)-\delta_i),$$
where $\delta_i=0$ if $G^f_i$ is equal to its ring and $\delta_i=2$ otherwise.
Furthermore, Lemma~\ref{lemma-obstacle} guarantees that $\el(f)\le 2$ and that if $\el(f)>0$, then $S_f$ is non-trivial.
Hence, summing these inequalities, we can strengthen (\ref{eq-partsfprel}) to
\begin{equation}\label{eq-partsf}
\sum_{i=1}^{k_f} w(G^f_i,\RR^f_i)\le w(f).
\end{equation}

As the sum of the elasticities of the faces of $G'$ is guaranteed to be at most $4$ by Lemma~\ref{lemma-obstacle},
inequalities (\ref{eq-cover}) and (\ref{eq-partsfprel}) give
\begin{align*}
w(G,\RR)&\le \sum_{f\in F(G')} (w(f)+\el(f)+\eta\cdot\surf(f))\label{eq-main}\\
&\le w(G',\RR)+4+\eta\sum_{f\in F(G')} \surf(f).\nonumber
\end{align*}
If $G'$ has a face that is neither open $2$-cell nor omnipresent, then
(\ref{cl-non2nonomni}) applied to $G_0$ and (\ref{ineq-flip}) imply that
$$w(G',\RR)\le \ell(\RR) + \eta\Bigl(\surf(g, |\RR|, t_0,t_1) - 7 - \sum_{f\in F(G')} \surf(f)\Bigr),$$
and consequently $G$ satisfies the outcome of Theorem~\ref{thm-maingen}.
Therefore, we can assume that all faces of $G'$ are either open $2$-cell or omnipresent.
Similarly, using (\ref{cl-rep1}) we can assume that if no face of $G'$ is omnipresent, then
all of them are semi-closed $2$-cell.

Suppose first that $G'$ has no omnipresent face.  Using (\ref{eq-cover}), (\ref{eq-partsf}), (\ref{ineq-flip}), and
applying Theorem~\ref{thm-maingen} inductively for $G_0$, we have
\begin{align*}
w(G,\RR)&=\sum_{f\in F(G')} \sum_{i=1}^{k_f} w(G^f_i,\RR^f_i)\\
&\le \sum_{f\in F(G')} w(f)\\
&=w(G',\RR)\le w(G_0,\RR_0)\\
&\le\ell(\RR)+\eta\cdot\surf(g,|\RR|,t_0,t_1),
\end{align*}
hence $G$ satisfies the outcome of Theorem~\ref{thm-maingen}.

It remains to consider the case that $G'$ has an omnipresent face $h$.
Then, every component of $G$ is a plane graph with one ring, and by
Lemma~\ref{lem:crit3conn}, we conclude that every face
of $G$ different from $h$ is closed $2$-cell.  Lemma~\ref{lemma-obstacle} guarantees $\el(h)\le 2$.
By (\ref{eq-cover}), (\ref{eq-partsfprel}), (\ref{eq-partsf}) and (\ref{ineq-flip}), we have
\begin{align*}
w(G,\RR)&\le w(h)+\el(h)+\eta\cdot\surf(h)+\sum_{f\in F(G'),f\neq h} w(f)\\
&=w(G',\RR)+\el(h)+\eta\cdot\surf(h)\le w(G_0,\RR_0)+2+\eta\cdot\surf(h).
\end{align*}
However, by (\ref{cl-omni}), we have $w(G_0,\RR_0)\le \ell(\RR_0)-2=\ell(\RR)-2+\eta(\surf(g,|\RR|,t_0,t_1)-\surf(h))$,
and the claim of the theorem follows.
\end{proof}

For use in a future paper of this series, let us formulate a corollary of Theorem~\ref{thm-maingen}.
Let $G$ be a graph embedded normally in a surface $\Sigma$.
For a real number $\eta$ and a face $f$ of $G$, let $w_\eta(f)=w(f)+\eta\cdot\surf(f)$.
We define $w_\eta(G)$ as the sum of $w_\eta(f)$ over the faces $f$ of $G$.

\begin{corollary}
There exists a constant $\eta>0$ such that the following holds.
Let $G$ be a triangle-free graph embedded in a surface $\Sigma$ without non-contractible $4$-cycles and with rings $\RR$.
Let $B$ be the union of the rings $\RR$.
If $G$ is $\RR$-critical, then $w_\eta(G)\le w_\eta(B)$.
\end{corollary}
\begin{proof}
Choose $\eta$ as the constant of Theorem~\ref{thm-maingen}.  Suppose first that $G$ is embedded in the disk with one ring of length $l\ge 4$.
Since $G$ is $\RR$-critical, it is connected and thus all faces of $G$ are open $2$-cell, and $w_\eta(f)=w(f)$ for every face $f$ of $G$.
Note furthermore that $w_\eta(B)=s(l)$.  Hence, $w_\eta(G)=w(G,\RR)\le s(l-2)\le w_\eta(B)$ by Corollary~\ref{cor-diskcase}.

Hence, we can assume that $\Sigma$ is not the disk, and thus $$w_\eta(B)=\ell(\RR)+\eta\cdot\surf(g(\Sigma),|\RR|, t_0(\RR), t_1(\RR)).$$
If all faces of the embedding of $G$ in $\Sigma$ are open $2$-cell, then $w_\eta(G)=w(G,\RR)$ and the claim follows from Theorem~\ref{thm-maingen}.
Otherwise, the claim follows from (\ref{cl-non2nonomni}) or (\ref{cl-omni}).
\end{proof}

\section{Forcing $3$-colorability by removing $O(g)$ vertices}\label{sec-linbound}

Finally, we show that in a triangle-free graph of genus $g$, it suffices to remove $O(g)$ vertices to make it $3$-colorable.

\begin{proof}[Proof of Theorem~\ref{thm-linbound}]
Let $\kappa$ be the constant of Theorem~\ref{thm:main}.  We let $\beta=\max(5\kappa,4)$.
We prove Theorem~\ref{thm-linbound} by induction on the number of vertices.  For $g=0$, the claim holds by Gr\"otzsch's theorem, hence assume that $g>0$.
Let $G$ be a triangle-free graph embedded in a surface $\Sigma$ of Euler genus $g$, such that Theorem~\ref{thm-linbound} holds for all graphs with less than $|V(G)|$
vertices.  Without loss of generality, we can assume that all the faces in the embedding of $G$ are open $2$-cell.

Suppose that $G$ contains a non-facial $4$-cycle $K$.
If $K$ is non-contractible, then $G-V(K)$ can be embedded in a surface of Euler genus at most $g-1$.  By the induction hypothesis,
there exists a set $X'$ of size at most $\beta(g-1)$ such that $G-V(K)-X'$ is $3$-colorable, and we can set $X=V(K)\cup X'$.

Hence, assume that $K$ is contractible.  Let $G'$ be obtained from $G$ by removing vertices and edges contained in the open disk
bounded by $K$.  By the induction hypothesis, there exists a set $X$ of size at most $\beta g$ such that $G'-X$ is $3$-colorable.
We claim that $G-X$ is also $3$-colorable.  Indeed, consider any $3$-coloring $\varphi$ of $G'-X$ and let $H$ be the subgraph of $G$ drawn in the closed disk bounded by $K$.
Let $\psi$ be a $3$-coloring of $K$ that matches $\varphi$ on $V(K)\setminus X$.  By Theorem~\ref{thm:six}, $\psi$ extends to a $3$-coloring $\varphi'$ of $H$.
Hence, $\varphi$ together with the restriction of $\varphi'$ to $V(H)\setminus X$ gives a $3$-coloring of $G-X$.

Therefore, we can assume that every $4$-cycle in $G$ bounds a face.  If $G$ is $4$-critical, then let $X=\{v\}$ for an arbitrary vertex $v\in V(G)$.  By the definition of a $4$-critical
graph, $G-X$ is $3$-colorable.  Hence, assume that $G$ is not a $4$-critical graph.

Let $G_0\subseteq G$ be a maximal subgraph of $G$ such that every component of $G_0$ is $4$-critical.  Let $X$ consist of all vertices of $G_0$ that are incident with a face of length greater than $4$.
Firstly, we bound the size of $X$.  Let $H_1$, \ldots, $H_k$ be the components of $G_0$.  For $i=1,\ldots, k$, let us obtain an embedding of $H_i$ in a surface $\Sigma_i$ such that every face
is open $2$-cell from the embedding of $H_i$ in $\Sigma$ as follows: as long as there exists a face $f$ that is not open $2$-cell, cut along a closed non-contractible curve contained in $f$
and cap the resulting holes by disks.  Since $H_i$ is $4$-critical, triangle-free and every $4$-cycle in $H_i$ bounds a face, Theorem~\ref{thm:main} implies that
$|X\cap V(H_i)|\le 5\kappa g(\Sigma_i)$.  Observe that since $H_1$, \ldots, $H_k$ are vertex-disjoint, we have $g\ge g(\Sigma_1)+\ldots+g(\Sigma_k)$, and thus $|X|\le 5\kappa g$ as required.

Consider any $i\in\{1,\ldots,k\}$.  If $V(H_i)\cap X=\emptyset$, then all faces of $H_i$ have length $4$.  Since all $4$-cycles in $G$ bound faces, it follows that $G=H_i$.  However, this contradicts the assumption that $G$ is not $4$-critical.
Therefore, $X$ intersects every component of $G_0$.

We claim that $G-X$ is $3$-colorable.  Indeed, if that is not the case, then $G-X$ contains a $4$-critical subgraph $G'$.  By the maximality of $G_0$, it follows that there exists $i\in\{1,\ldots, k\}$
such that $G'$ intersects $H_i$.  Since $V(H_i)\cap X\neq\emptyset$, we have $G'\neq H_i$, and since $G'$ and $H_i$ are $4$-critical, it follows that $G'\not\subseteq H_i$.
Since $G'$ is $4$-critical, it is connected, and thus there exists $uv\in E(G)$ such that $u\in V(G')\setminus V(H_i)$ and $v\in V(G')\cap V(H_i)$.
Let $f$ be the face of $H_i$ containing $u$, and note that $v$ is incident with $f$.  Since $v\in V(G')$, we have $v\not\in X$, and thus $f$ is a $4$-face.  However,
then the $4$-cycle bounding the face $f$ separates $u$ from all vertices of $X$ in $G$, which contradicts the assumption that all $4$-cycles in $G$ bound faces.
\end{proof}

\bibliographystyle{acm}
\bibliography{ctvrty}

\begin{thebibliography}{10}

\bibitem{aksenov}
{\sc Aksionov, V.~A.}
\newblock On continuation of $3$-colouring of planar graphs.
\newblock {\em Diskret. Anal. Novosibirsk 26\/} (1974), 3--19.
\newblock In Russian.

\bibitem{DvoKawTho}
{\sc Dvo\v{r}\'ak, Z., Kawarabayashi, K., and Thomas, R.}
\newblock Three-coloring triangle-free planar graphs in linear time.
\newblock {\em Trans. on Algorithms 7\/} (2011), article no. 41.

\bibitem{dkt}
{\sc Dvo\v{r}\'ak, Z., Kr\'al', D., and Thomas, R.}
\newblock Coloring planar graphs with triangles far apart.
\newblock {\em ArXiv 0911.0885v1\/} (Nov. 2009).

\bibitem{trfree5}
{\sc Dvo\v{r}\'ak, Z., Kr\'al', D., and Thomas, R.}
\newblock Three-coloring triangle-free graphs on surfaces {V}. {C}oloring
  planar graphs with distant anomalies.
\newblock {\em ArXiv 0911.0885v3\/} (Jan. 2016).

\bibitem{trfree2}
{\sc Dvo\v{r}\'ak, Z., Kr\'al', D., and Thomas, R.}
\newblock Three-coloring triangle-free graphs on surfaces {II}. $4$-critical
  graphs in a disk.
\newblock {\em Journal of Combinatorial Theory, Series B 132\/} (2018), 1--46.

\bibitem{trfree3}
{\sc Dvo\v{r}\'ak, Z., Kr\'al', D., and Thomas, R.}
\newblock Three-coloring triangle-free graphs on surfaces {III}. {G}raphs of
  girth five.
\newblock {\em Journal of Combinatorial Theory, Series B 145\/} (2020),
  376--432.

\bibitem{gimbel}
{\sc Gimbel, J., and Thomassen, C.}
\newblock Coloring graphs with fixed genus and girth.
\newblock {\em Trans. Amer. Math. Soc. 349\/} (1997), 4555--4564.

\bibitem{grotzsch1959}
{\sc Gr{\"o}tzsch, H.}
\newblock Ein {D}reifarbensatz f\"{u}r dreikreisfreie {N}etze auf der {K}ugel.
\newblock {\em Math.-Natur. Reihe 8\/} (1959), 109--120.

\bibitem{conj-havel}
{\sc Havel, I.}
\newblock On a conjecture of {B}. {G}r\"unbaum.
\newblock {\em J. Combin. Theory, Ser.~B 7\/} (1969), 184--186.

\bibitem{locplanq}
{\sc Hutchinson, J.~P.}
\newblock Three-coloring graphs embedded on surfaces with all faces even-sided.
\newblock {\em J. Combin. Theory, Ser.~B 65\/} (1995), 139--155.

\bibitem{kawthorem}
{\sc Kawarabayashi, K., and Thomassen, C.}
\newblock From the plane to higher surfaces.
\newblock {\em Journal of Combinatorial Theory, Series B 102}, 4 (2012),
  852--868.

\bibitem{koyan}
{\sc Kostochka, A.~V., and Yancey, M.}
\newblock Ore's conjecture for $k=4$ and gr{\"{o}}tzsch's theorem.
\newblock {\em Combinatorica 34}, 3 (2014), 323--329.

\bibitem{rs7}
{\sc Robertson, N., and Seymour, P.~D.}
\newblock Graph {M}inors {VII}. {D}isjoint paths on a surface.
\newblock {\em J. Combin. Theory, Ser.~B 45\/} (1988), 212--254.

\bibitem{tw-klein}
{\sc Thomas, R., and Walls, B.}
\newblock Three-coloring {K}lein bottle graphs of girth five.
\newblock {\em J. Combin. Theory, Ser.~B 92\/} (2004), 115--135.

\bibitem{thom-torus}
{\sc Thomassen, C.}
\newblock Gr{\"o}tzsch's 3-color theorem and its counterparts for the torus and
  the projective plane.
\newblock {\em J. Combin. Theory, Ser.~B 62\/} (1994), 268--279.

\bibitem{thomassen1995-34}
{\sc Thomassen, C.}
\newblock 3-list-coloring planar graphs of girth 5.
\newblock {\em J. Combin. Theory, Ser.~B 64\/} (1995), 101--107.

\bibitem{thomassen-surf}
{\sc Thomassen, C.}
\newblock The chromatic number of a graph of girth 5 on a fixed surface.
\newblock {\em J. Combin. Theory, Ser.~B 87\/} (2003), 38--71.

\bibitem{ThoShortlist}
{\sc Thomassen, C.}
\newblock A short list color proof of {G}r{\"o}tzsch's theorem.
\newblock {\em J. Combin. Theory, Ser.~B 88\/} (2003), 189--192.

\end{thebibliography}

\end{document}